# AN ADAPTATION THEORY FOR NONPARAMETRIC CONFIDENCE INTERVALS[1]

By T. Tony Cai and Mark G. Low

*University of Pennsylvania*

A nonparametric adaptation theory is developed for the construction of confidence intervals for linear functionals. A between class modulus of continuity captures the expected length of adaptive confidence intervals. Sharp lower bounds are given for the expected length and an ordered modulus of continuity is used to construct adaptive confidence procedures which are within a constant factor of the lower bounds. In addition, minimax theory over nonconvex parameter spaces is developed.

**1. Introduction.** The problem of estimating a linear functional occupies a central position in nonparametric function estimation. It is most complete in the Gaussian settings:

$$(1) \qquad dY(t) = f(t)\,dt + n^{-1/2}\,dW(t), \qquad -\tfrac{1}{2} \le t \le \tfrac{1}{2},$$

where $W(t)$ is standard Brownian motion and

$$(2) \qquad Y(i) = f(i) + n^{-1/2} z_i, \qquad i \in \mathcal{M},$$

where $z_i$ are i.i.d. standard normal random variables and $\mathcal{M}$ is a finite or countably infinite index set. In particular, minimax estimation theory has been well developed in Ibragimov and Hasminskii (1984), Donoho and Liu (1991) and Donoho (1994).

Confidence sets also play a fundamental role in statistical inference. In the context of nonparametric function estimation variable size confidence intervals, bands and balls have received particular attention recently. For any confidence set there are two main interrelated issues which need to

Received February 2003; revised December 2003.
[1]Supported in part by NSF Grant DMS-03-06576.

be considered together, coverage probability and the expected size of the confidence set.

One common technique for constructing confidence bands and intervals is through the bootstrap. In this context it has been noted that intervals based on the bootstrap often have poor coverage probability. See, for example, Hall (1992) and Härdle and Marron (1991). Picard and Tribouley (2000) construct adaptive confidence intervals for functions at a point using a wavelet method which achieve optimal coverage accuracy up to a logarithmic factor although in this case the issue of optimal expected length is not addressed. On the other hand Li (1989), Beran and Dümbgen (1998) and Genovese and Wasserman (2002) have constructed confidence balls which guarantee coverage probability. Closer to the present work, adaptive confidence bands have been constructed in the special case of shape restricted functions. In this context Hengartner and Stark (1995) and Dümbgen (1998) give a variable width confidence band which adapts to local smoothness while maintaining a given level of coverage probability.

In this paper we focus on the construction of confidence intervals for linear functionals which adapt to the unknown function. This adaptation problem can be made precise by considering collections of parameter spaces $\{\mathcal{F}_j, j \in J\}$, where $J$ is some index set. For such a collection of parameter spaces the confidence interval should have a given coverage probability over the union of the parameter spaces. Subject to this constraint the goal is to minimize the maximum expected length simultaneously over each of the parameter spaces.

For example, consider the simple and most easily explained case of two nested spaces, $\mathcal{F}_1 \subseteq \mathcal{F}$. An adaptive confidence interval must attain optimal expected length performance over both $\mathcal{F}_1$ and $\mathcal{F}$ while satisfying a given coverage probability over $\mathcal{F}$. More specifically write $\mathcal{I}_{\alpha,\mathcal{F}}$ for the collection of all confidence intervals which cover the linear functional $Tf$ with minimum coverage probability of at least $1 - \alpha$ over the parameter space $\mathcal{F}$. Denote by $L(CI, \mathcal{G}) = \sup_{f \in \mathcal{G}} E_f(L(CI))$ the maximum expected length of a confidence interval $CI$ over $\mathcal{G}$ where $L(CI)$ is the length of the $CI$. Then a benchmark for the evaluation of the maximum expected length over $\mathcal{F}_1$ for any $CI \in \mathcal{I}_{\alpha,\mathcal{F}}$ is given by

(3) $$L_\alpha^*(\mathcal{F}_1, \mathcal{F}) = \inf_{CI \in \mathcal{I}_{\alpha,\mathcal{F}}} L(CI, \mathcal{F}_1).$$

In particular, when $\mathcal{F}_1 = \mathcal{F}$ set $L_\alpha^*(\mathcal{F}) = L_\alpha^*(\mathcal{F}, \mathcal{F})$, which gives the minimax expected length of confidence intervals of level $1 - \alpha$ over $\mathcal{F}$. For convex $\mathcal{F}$, Donoho (1994) constructed fixed length intervals centered at affine estimators which have length within a small constant factor of $L_\alpha^*(\mathcal{F})$.

The major result in the present paper is the construction of confidence intervals which have expected length within a constant factor of $L_\alpha^*(\mathcal{F}_j, \mathcal{F})$ simultaneously over a collection of convex parameter spaces $\mathcal{F}_j$ where $\mathcal{F} = \cup \mathcal{F}_j$. The construction of such intervals is general and is applicable to collections of arbitrary convex parameter spaces. It is shown in Cai and Low (2003) that in particular cases, such as collections of convex functions, the general procedure can be modified to yield simple and easily implementable procedures.

The main technical tools used in the derivation of the general adaptive confidence intervals are geometric quantities, the ordered and between class moduli of continuity which are defined as follows. For a linear functional $T$ and parameter spaces $\mathcal{F}$ and $\mathcal{G}$ there are ordered moduli of continuity $\omega(\varepsilon, \mathcal{F}, \mathcal{G})$ associated with the Gaussian models (1) and (2) defined by

(4) $\qquad \omega(\varepsilon, \mathcal{F}, \mathcal{G}) = \sup\{Tg - Tf : \|g - f\|_2 \leq \varepsilon; f \in \mathcal{F}, g \in \mathcal{G}\},$

where $\|\cdot\|_2$ is the $L_2(-\frac{1}{2}, \frac{1}{2})$ function norm in the white noise model (1) and the $\ell_2$ sequence norm over the index set $\mathcal{M}$ in the Gaussian model (2). As we shall give a unified treatment of both models it is convenient in the notation used throughout the paper not to distinguish the function norm and the sequence norm. It is implicit that for results concerning the white noise model (1) the notation $\|\cdot\|_2$ always refers to the $L_2$ function norm whereas for the sequence model (2) it always refers to the $\ell_2$ sequence norm. When $\mathcal{G} = \mathcal{F}$, $\omega(\varepsilon, \mathcal{F}, \mathcal{F})$ is the modulus of continuity over $\mathcal{F}$ introduced by Donoho and Liu (1991) and will be denoted by $\omega(\varepsilon, \mathcal{F})$.

For two parameter spaces $\mathcal{F}$ and $\mathcal{G}$ and a given linear functional $T$, the between class modulus of continuity is defined as $\omega_+(\varepsilon, \mathcal{F}, \mathcal{G}) = \max\{\omega(\varepsilon, \mathcal{F}, \mathcal{G}), \omega(\varepsilon, \mathcal{G}, \mathcal{F})\}$, or equivalently

(5) $\qquad \omega_+(\varepsilon, \mathcal{F}, \mathcal{G}) = \sup\{|Tg - Tf| : \|g - f\|_2 \leq \varepsilon; f \in \mathcal{F}, g \in \mathcal{G}\}.$

The between class and ordered moduli were first introduced in Cai and Low (2002) in the context of adaptive estimation under mean squared error where they were shown to be instrumental in characterizing the possible degree of adaptability over two convex classes $\mathcal{F}$ and $\mathcal{G}$ in the same way that the modulus of continuity $\omega(\varepsilon, \mathcal{F})$ used by Donoho and Liu (1991) and Donoho (1994) captures the minimax difficulty of estimation over a single convex parameter space $\mathcal{F}$.

The paper is organized as follows. Section 2 covers adaptation over two convex parameter spaces $\mathcal{F}_1$ and $\mathcal{F}_2$ where the theory is most easily understood. A lower bound based on the between class modulus as defined in (5) is given for $L_\alpha^*(\mathcal{F}_1, \mathcal{F})$ where $\mathcal{F} = \mathcal{F}_1 \cup \mathcal{F}_2$. An adaptive confidence interval attaining this bound is also constructed by using the ordered moduli as given in (4). Various examples are used to illustrate the adaptation theory.

More generally let $\{\mathcal{F}_j, j \in J\}$ be a collection of convex parameter spaces with nonempty intersections and let $\mathcal{F} = \cup \mathcal{F}_j$. The goal is then to simultaneously minimize $L(CI, \mathcal{F}_j)$ for confidence intervals $CI \in \mathcal{I}_{\alpha, \mathcal{F}}$. For each parameter space $\mathcal{F}_j$, $L^*_\alpha(\mathcal{F}_j, \mathcal{F})$ provides a lower bound on the maximum expected length over $\mathcal{F}_j$ for any $CI \in \mathcal{I}_{\alpha, \mathcal{F}}$. In Section 3 a complete treatment is given for nested $\mathcal{F}_j$, possibly infinite in number. For any collection of nested convex parameter spaces a variable length confidence interval is constructed which for a given level of coverage has expected length within a constant factor of the minimum expected length simultaneously over all parameter spaces in the collection.

Section 4 treats the case of a general finite collection of convex parameter spaces. A more complicated procedure results in an interval which also has expected length within a constant factor of the minimum expected length although the constant factor now depends on the number of parameter spaces in the collection. Finally in Section 5 it is shown, by example, that the rate of growth in this constant factor as a function of the number of parameter spaces cannot in general be avoided. In addition, the adaptation theory developed in this paper is used to extend the minimax theory to a finite union of convex parameter spaces. This extension is given in Section 5.

**2. Adaptation over two parameter spaces.** In this section we consider adaptation over two parameter spaces. For the development of this theory, it is convenient for a given $\alpha$ to provide a benchmark for the maximum expected length over $\mathcal{F}_1$ of confidence intervals with a given coverage probability of $1 - \alpha$ over $\mathcal{F} = \mathcal{F}_1 \cup \mathcal{F}_2$, namely to provide a lower bound for $L^*_\alpha(\mathcal{F}_1, \mathcal{F})$ as defined in (3). This benchmark is given in Section 2.1 for arbitrary parameter spaces.

We give a complete treatment of adaptation when the two parameter spaces are convex. In this case adaptive intervals attaining the lower bound given in Section 2.1 are constructed. The adaptive procedure is given in Section 2.2. Examples illustrating the theory are given in Section 2.3.

It is convenient to write $a_l \asymp b_l$ whenever $0 < \liminf a_l/b_l \leq \limsup a_l/b_l < \infty$, where $l$ ranges over either a continuous or discrete index set.

2.1. *Lower bound on the length of confidence intervals.* The following simple two-point Normal mean problem is the basis for a surprisingly useful general lower bound on the expected length of $1 - \alpha$ level confidence intervals. We shall see later that the two-point bound is easy to apply for adaptation theory because each point can be chosen to lie in different parameter spaces. Previous work on confidence intervals for bounded Normal means as in Pratt (1961), Zeytinoglu and Mintz (1984) and Stark (1992) is useful for minimax theory but it is not applicable for general adaptation problems.

Let $X \sim N(\theta, \sigma^2)$ and suppose that $\theta \in \Theta = \{\theta_0, \theta_1\}$ where $\theta_0 < \theta_1$. Consider the following simple statistical decision theory problem: construct confidence intervals $CI(X)$ for $\theta$ which have smallest expected length under $\theta_0$ subject to the coverage constraint

$$P_\theta(\theta \in CI(X)) \geq 1 - \alpha \quad \text{for } \theta \in \Theta.$$

Throughout the paper set $z_\alpha = \Phi^{-1}(1-\alpha)$ where $\Phi$ is the cumulative density function of a standard Normal distribution. In addition write $L(CI)$ for the length of a confidence interval $CI$.

PROPOSITION 1. *Let $X \sim N(\theta, \sigma^2)$ and suppose that $\theta \in \Theta = \{\theta_0, \theta_1\}$ where $\theta_0 < \theta_1$. Let $CI(X)$ be a $1 - \alpha$ level confidence interval for $\theta$. Then*

$$(6) \quad E_{\theta_i} L(CI(X)) \geq (\theta_1 - \theta_0)\left(1 - \alpha - \Phi\left(\frac{\theta_1 - \theta_0}{\sigma} - z_\alpha\right)\right)_+$$

*for $i = 0, 1$. Moreover there exists a confidence interval which attains the lower bounds simultaneously for both $i = 0$ and $i = 1$.*

PROOF. It is clear that it suffices to consider confidence intervals $CI(X)$ of three possible forms: $[\theta_0, \theta_1]$, $\{\theta_0\}$ and $\{\theta_1\}$. The problem is then to minimize $P_{\theta_0}(CI(X) = [\theta_0, \theta_1])$ subject to the constraints $P_{\theta_0}(CI(X) = \{\theta_1\}) \leq \alpha$ and $P_{\theta_1}(CI(X) = \{\theta_0\}) \leq \alpha$.

It follows from the Neyman–Pearson lemma that, subject to the constraint that $P_{\theta_1}(CI(X) = \{\theta_0\}) \leq \alpha$,

$$P_{\theta_0}(CI(X) = \{\theta_0\}) \leq \Phi\left(\frac{\theta_1 - \theta_0}{\sigma} - z_\alpha\right).$$

Hence

$$E_{\theta_0} L(CI(X)) = (\theta_1 - \theta_0) P_{\theta_0}(CI(X) = [\theta_0, \theta_1])$$
$$= (\theta_1 - \theta_0)(1 - P_{\theta_0}(CI(X) = \{\theta_1\}) - P_{\theta_0}(CI(X) = \{\theta_0\}))$$
$$\geq (\theta_1 - \theta_0)\left(1 - \alpha - \Phi\left(\frac{\theta_1 - \theta_0}{\sigma} - z_\alpha\right)\right).$$

The bound for $\theta_1$ follows similarly.

It is easy to see that an interval attaining the lower bound for $\theta_0$ and $\theta_1$ is given by

$$CI(X) = \begin{cases} \{\theta_0\}, & \text{if } X \leq \theta_1 - z_\alpha \sigma, \\ [\theta_0, \theta_1], & \text{if } \theta_1 - z_\alpha \sigma < X < \theta_0 + z_\alpha \sigma, \\ \{\theta_1\}, & \text{if } X \geq \theta_0 + z_\alpha \sigma, \end{cases}$$

when $\theta_1 - z_\alpha \sigma < \theta_0 + z_\alpha \sigma$. Otherwise set

$$CI(X) = \begin{cases} \{\theta_0\}, & \text{if } X \leq \frac{\theta_0 + \theta_1}{2}, \\ \{\theta_1\}, & \text{if } X > \frac{\theta_0 + \theta_1}{2}. \end{cases}$$

In this case the confidence interval always has zero length and coverage of at least $1 - \alpha$. $\square$

Based on the two-point bound given in Proposition 1 the following theorem gives a lower bound for infinite-dimensional Gaussian models.

THEOREM 1. *Let $0 < \alpha < \frac{1}{2}$ and let $\mathcal{F}_1 \subseteq \mathcal{F}$ be two parameter spaces. Then*

$$L_\alpha^*(\mathcal{F}_1, \mathcal{F}) \geq \left(\frac{1}{2} - \alpha\right) \omega_+ \left(\frac{z_\alpha}{\sqrt{n}}, \mathcal{F}_1, \mathcal{F}\right) \tag{7}$$

*where $L_\alpha^*(\mathcal{F}_1, \mathcal{F})$ is defined in* (3) *and $\omega_+(\varepsilon, \mathcal{F}_1, \mathcal{F})$ is the between class modulus as given in* (5).

PROOF. We shall focus on the proof for the white noise with drift model (1). The proof for the sequence model (2) is analogous. Fix $\varepsilon > 0$. For any $\delta > 0$ there are functions $f_1 \in \mathcal{F}_1$ and $f_2 \in \mathcal{F}$ such that

$$|Tf_2 - Tf_1| \geq \omega_+\left(\frac{\varepsilon}{\sqrt{n}}, \mathcal{F}_1, \mathcal{F}\right) - \delta$$

and such that

$$\|f_2 - f_1\|_2 \leq \frac{\varepsilon}{\sqrt{n}}.$$

Denote by $P_i$ the probability measure associated with the white noise process

$$dY(t) = f_i(t)\,dt + \frac{1}{\sqrt{n}}\,dW(t), \qquad -\tfrac{1}{2} \leq t \leq \tfrac{1}{2}, i = 1, 2.$$

Let $\beta_n = n\|f_1 - f_2\|_2^2$. Then a sufficient statistic for the family of measures $\{P_i : i = 1, 2\}$ is given by the log-likelihood ratio $S_n = \log(dP_2/dP_1)$ with

$$S_n \sim \begin{cases} N-\left(\dfrac{\beta_n}{2}, \beta_n\right) & \text{under } P_1, \\ N\left(\dfrac{\beta_n}{2}, \beta_n\right) & \text{under } P_2. \end{cases}$$

An equivalent sufficient statistic is thus given by

$$Q_n = \frac{Tf_1 + Tf_2}{2} + \frac{Tf_2 - Tf_1}{\beta_n} \cdot S_n$$

where

$$Q_n \sim \begin{cases} N\left(Tf_1, \dfrac{(Tf_2 - Tf_1)^2}{\beta_n}\right) & \text{under } P_1, \\ N\left(Tf_2, \dfrac{(Tf_2 - Tf_1)^2}{\beta_n}\right) & \text{under } P_2. \end{cases}$$

It follows from Proposition 1 that for any confidence interval $CI(Q_n)$ based on $Q_n$,

$$E_{f_1} L(CI(Q_n)) \geq |Tf_2 - Tf_1|\left(1 - \alpha - \Phi\left(\frac{|Tf_2 - Tf_1|}{\sigma} - z_\alpha\right)\right)_+$$

where $\sigma = \frac{|Tf_2 - Tf_1|}{\sqrt{\beta_n}}$. Hence

$$E_{f_1} L(CI(Q_n)) \geq |Tf_2 - Tf_1|(1 - \alpha - \Phi(\sqrt{n}\|f_2 - f_1\|_2 - z_\alpha))_+$$
$$\geq \left(\omega_+\left(\frac{\varepsilon}{\sqrt{n}}, \mathcal{F}_1, \mathcal{F}\right) - \delta\right)(1 - \alpha - \Phi(\varepsilon - z_\alpha))_+.$$

Letting $\delta \to 0$, it follows that for any $\varepsilon > 0$,

$$L(CI(Q_n), \mathcal{F}_1) \geq \omega_+\left(\frac{\varepsilon}{\sqrt{n}}, \mathcal{F}_1, \mathcal{F}\right)(1 - \alpha - \Phi(\varepsilon - z_\alpha))_+.$$

By the sufficiency of $Q_n$, it follows that for any confidence interval $CI \in \mathcal{I}_{\alpha, \mathcal{F}}$

(8) $$L(CI, \mathcal{F}_1) \geq \sup_{\varepsilon > 0} \omega_+\left(\frac{\varepsilon}{\sqrt{n}}, \mathcal{F}_1, \mathcal{F}\right)(1 - \alpha - \Phi(\varepsilon - z_\alpha))_+.$$

The theorem follows on taking $\varepsilon = z_\alpha$. $\square$

REMARK 1. Although the primary use of this theorem is for adaptive confidence intervals, it can also be used to show that from a minimax point of view there is relatively little to gain by using variable length intervals. In the minimax setting Donoho (1994) showed that over a given convex parameter space $\mathcal{F}$, fixed length confidence intervals for a linear functional $Tf$ with coverage of at least $1 - \alpha$ must have maximum length at least $2\omega(\frac{2z_\alpha}{\sqrt{n}}, \mathcal{F})$ and that fixed length confidence intervals can be centered on affine estimators with maximum length at most $2\omega(\frac{2z_{\alpha/2}}{\sqrt{n}}, \mathcal{F})$. By taking $\mathcal{F}_1 = \mathcal{F}$, Theorem 1 yields that the minimax expected length of a $1 - \alpha$ level confidence interval over any parameter space $\mathcal{F}$ satisfies

(9) $$L_\alpha^*(\mathcal{F}) \geq \left(\frac{1}{2} - \alpha\right)\omega\left(\frac{z_\alpha}{\sqrt{n}}, \mathcal{F}\right).$$

This shows that for any given $\alpha < 1/2$ the optimal variable length confidence intervals must have maximum expected length at least a fixed constant factor of the length of the shortest fixed length confidence interval when the parameter space $\mathcal{F}$ is convex.

2.2. *Adaptive confidence interval.* There are at least two natural ways to define adaptive confidence intervals over a collection of convex parameter spaces $\{\mathcal{F}_i, i = 1, \ldots, k\}$. Let $\mathcal{F} = \bigcup_{i=1}^k \mathcal{F}_i$. Call a confidence interval $CI \in \mathcal{I}_{\alpha,\mathcal{F}}$ *adaptive* over the collection $\{\mathcal{F}_i, i = 1, \ldots, k\}$ if, for all $1 \leq i \leq k$,

$$L(CI, \mathcal{F}_i) \leq C_i(\alpha)\omega_+\left(\frac{z_\alpha}{\sqrt{n}}, \mathcal{F}_i, \mathcal{F}\right), \tag{10}$$

where $C_i(\alpha)$ are constants depending on $\alpha$ only. In other words a confidence interval which adapts over the parameter spaces $\mathcal{F}_i$ attains the lower bound given in Theorem 1 for each $i$ while maintaining coverage over $\mathcal{F}$. We shall show that such adaptive confidence intervals can always be constructed when $k$ is finite.

It is also reasonable, in light of the minimax discussion given above, to term a confidence interval $CI \in \mathcal{I}_{\alpha,\mathcal{F}}$ adaptive over the collection of parameter spaces $\mathcal{F}_i$ if, for all $1 \leq i \leq k$,

$$L(CI, \mathcal{F}_i) \leq C_i(\alpha)\omega\left(\frac{z_\alpha}{\sqrt{n}}, \mathcal{F}_i\right) \tag{11}$$

where $C_i(\alpha)$ are constants depending on $\alpha$ only. We shall call such a confidence interval *strongly adaptive*. It is clear that a confidence interval which is strongly adaptive is also adaptive. However strongly adaptive confidence intervals do not always exist. Low (1997) has given examples where $L^*_\alpha(\mathcal{F}_1, \mathcal{F}) \gg L^*_\alpha(\mathcal{F}_1)$, in which case strongly adaptive estimators do not exist. Other examples are given in Section 2.3 and throughout the paper. On the other hand, when $L^*_\alpha(\mathcal{F}_1, \mathcal{F}) \asymp L^*_\alpha(\mathcal{F}_1)$ strongly adaptive estimators do exist and any estimator which is adaptive is also strongly adaptive.

In this section the focus is on adaptation over two parameter spaces where the theory is most easily understood. For two parameter spaces $\mathcal{F}_1 \subseteq \mathcal{F}$, Theorem 1 gives a lower bound for the maximum expected length over $\mathcal{F}_1$ of confidence intervals with guaranteed coverage over $\mathcal{F}$. We now show that the lower bound can in fact be attained within a constant factor not depending on $n$ when $\mathcal{F}_1$ is convex and $\mathcal{F}$ is the union of $\mathcal{F}_1$ and another convex set $\mathcal{F}_2$.

Let $\{\mathcal{F}_1, \mathcal{F}_2\}$ be a pair of convex parameter spaces with nonempty intersection and let $\mathcal{F} = \mathcal{F}_1 \cup \mathcal{F}_2$. Our first objective is to construct a confidence interval for a linear functional $Tf$ which has guaranteed coverage probability of $1 - \alpha$ over $\mathcal{F}$ and has maximum expected length over $\mathcal{F}_1$ within a constant factor of the lower bound given in Theorem 1, namely, for any $CI \in \mathcal{I}_{\alpha,\mathcal{F}}$,

$$L(CI, \mathcal{F}_1) \geq \left(\frac{1}{2} - \alpha\right)\omega_+\left(\frac{z_\alpha}{\sqrt{n}}, \mathcal{F}_1, \mathcal{F}\right). \tag{12}$$

The construction of the adaptive confidence interval relies on the ordered modulus $\omega(\varepsilon, \mathcal{F}_i, \mathcal{F}_j)$ as given in (4). For $1 \leq i, j \leq 2$, set

$$\omega_{i,j} = \omega\left(\frac{z_{\alpha/2}}{\sqrt{n}}, \mathcal{F}_i, \mathcal{F}_j\right).$$

Cai and Low (2004) give an algorithm for the construction of a linear estimator $\hat{T}_{i,j}$ which has variance bounded by

(13) $$\operatorname{Var}(\hat{T}_{i,j}) \leq \frac{1}{z_{\alpha/2}^2} \omega_{i,j}^2$$

and bias which satisfies

(14) $$\inf_{f \in \mathcal{F}_j} (E(\hat{T}_{i,j}) - Tf) \geq -\tfrac{1}{2}\omega_{i,j}$$

and

(15) $$\sup_{f \in \mathcal{F}_i} (E(\hat{T}_{i,j}) - Tf) \leq \tfrac{1}{2}\omega_{i,j}.$$

We shall use the linear estimators $\hat{T}_{i,j}$ to construct a confidence interval which has guaranteed coverage probability over $\mathcal{F}$ and which also has expected length over $\mathcal{F}_1$ within a constant factor of the lower bound given by (26). For $j = 1$ and 2 define the confidence intervals $CI_{j,\alpha}^*$ by

(16) $$CI_{j,\alpha}^* = \left[\min_{i=1,2}\{\hat{T}_{i,j} - \tfrac{3}{2}\omega_{i,j}\}, \max_{i=1,2}\{\hat{T}_{j,i} + \tfrac{3}{2}\omega_{j,i}\}\right].$$

The following result shows that the confidence interval $CI_{1,\alpha}^*$ attains the lower bound on the maximum expected length over $\mathcal{F}_1$ given in (7) within a constant factor not depending on $n$ and satisfies the constraint that it has the minimum coverage of $1 - \alpha$ for all $f \in \mathcal{F}$.

LEMMA 1. *Let $\mathcal{F}_1$ and $\mathcal{F}_2$ be convex parameter spaces with $\mathcal{F}_1 \cap \mathcal{F}_2 \neq \varnothing$ and let $\mathcal{F} = \mathcal{F}_1 \cup \mathcal{F}_2$. Let the interval $CI_{j,\alpha}^*$ be defined as in (16) for $j = 1$ and 2. Then $CI_{j,\alpha}^* \in \mathcal{I}_{\alpha, \mathcal{F}}$ and $CI_{j,\alpha}^*$ has expected length over $\mathcal{F}_j$ which satisfies*

(17) $$L(CI_{j,\alpha}^*, \mathcal{F}_j) \leq \left\{\frac{9}{z_{\alpha/2}} + 4\right\} \omega_+\left(\frac{z_{\alpha/2}}{\sqrt{n}}, \mathcal{F}_j, \mathcal{F}\right).$$

Lemma 1 follows from the proof of Proposition 4 given in Section 4.1.

REMARK 2. Theorem 1 and Lemma 1 together show that under the conditions of Lemma 1,

(18) $$L_\alpha^*(\mathcal{F}_1, \mathcal{F}) \asymp \omega_+\left(\frac{z_{\alpha/2}}{\sqrt{n}}, \mathcal{F}_1, \mathcal{F}\right).$$

Although the interval $CI^*_{1,\alpha}$ has guaranteed coverage probability over $\mathcal{F}$ and optimal expected length over $\mathcal{F}_1$, it may not have optimal expected length over $\mathcal{F}$ because the expected length over $\mathcal{F}_2$ is not controlled. On the other hand, by symmetry $CI^*_{2,\alpha}$ has guaranteed coverage probability over $\mathcal{F}$ and optimal expected length over $\mathcal{F}_2$. By Bonferroni, the confidence interval $CI^*_\alpha = CI^*_{1,\alpha/2} \cap CI^*_{2,\alpha/2}$ also has coverage probability of at least $1 - \alpha$ and so $CI^*_\alpha \in \mathcal{I}_{\alpha,\mathcal{F}}$. Furthermore, it is easy to see that it has optimal expected length over both $\mathcal{F}_1$ and $\mathcal{F}_2$ and hence also over $\mathcal{F}$. In other words the confidence interval $CI^*_\alpha$ is a $1 - \alpha$ level adaptive confidence interval over $\mathcal{F}_1$ and $\mathcal{F}_2$.

PROPOSITION 2. *Let $\mathcal{F}_1$ and $\mathcal{F}_2$ be convex parameter spaces with $\mathcal{F}_1 \cap \mathcal{F}_2 \neq \varnothing$ and let $\mathcal{F} = \mathcal{F}_1 \cup \mathcal{F}_2$. Let the interval $CI^*_{j,\alpha}$ be defined as in (16) and let $CI^*_\alpha = CI^*_{1,\alpha/2} \cap CI^*_{2,\alpha/2}$. Then $CI^*_\alpha$ is a $1 - \alpha$ level adaptive confidence interval over $\mathcal{F}_1$ and $\mathcal{F}_2$. That is, $CI^*_\alpha \in \mathcal{I}_{\alpha,\mathcal{F}}$ and for both $j = 1$ and 2,*

$$(19) \qquad L^*_\alpha(\mathcal{F}_j, \mathcal{F}) \leq L(CI^*_\alpha, \mathcal{F}_j) \leq C(\alpha) L^*_\alpha(\mathcal{F}_j, \mathcal{F})$$

*where $C(\alpha)$ is a constant depending only on $\alpha$. Consequently $L(CI^*_\alpha, \mathcal{F}_j) \asymp \omega_+(\frac{z_{\alpha/2}}{\sqrt{n}}, \mathcal{F}_j, \mathcal{F})$.*

REMARK 3. It is shown in Cai and Low (2004) that the ordered modulus is concave. It follows that, if $b \geq 1$, then for all $\varepsilon > 0$,

$$(20) \qquad \begin{aligned} \omega_+(b\varepsilon, \mathcal{F}, \mathcal{G}) &= \max(\omega(b\varepsilon, \mathcal{F}, \mathcal{G}), \omega(b\varepsilon, \mathcal{G}, \mathcal{F})) \\ &\leq \max(b\omega(\varepsilon, \mathcal{F}, \mathcal{G}), b\omega(\varepsilon, \mathcal{G}, \mathcal{F})) \\ &\leq b\omega_+(\varepsilon, \mathcal{F}, \mathcal{G}). \end{aligned}$$

It then follows from the bounds given in (7) and (17) and inequality (20) that the constant $C(\alpha)$ in (19) can be taken as

$$C(\alpha) = \frac{9 + 4z_{\alpha/4}}{(1/2 - \alpha)z_\alpha}.$$

2.3. *Discussion.* In nonparametric function estimation the goal of adaptive estimation is often framed in terms of achieving optimality results simultaneously over a collection of parameter spaces $\{\mathcal{F}_j\}$. The benchmark for success is given by how well one could do if the parameter space is completely specified. We termed any such confidence interval strongly adaptive.

So far, attention has focused on constructing adaptive confidence procedures which attain the lower bound on expected length given in Theorem 1. This bound gives the best one can do in this adaptive confidence interval problem. The lower bound however may differ quite dramatically from the minimax expected length if the parameter space $\mathcal{F}_j$ is prespecified. In particular suppose, as is common, that the between class modulus of continuity is Hölderian. That is, the modulus satisfies

$$\omega_+(\varepsilon, \mathcal{F}_i, \mathcal{F}_j) = C_{i,j}\varepsilon^{q_{i,j}}(1+o(1)), \qquad 1 \leq i, j \leq 2,$$

for some constants $C_{i,j} > 0$ and $0 < q_{i,j} \leq 1$. Such is the case in the examples given in Section 3.2 and also in many other commonly treated problems. When the modulus $\omega_+(\varepsilon, \mathcal{F}, \mathcal{G})$ is Hölderian write $q(\mathcal{F}, \mathcal{G})$ for the exponent of the modulus. That is,

$$\omega_+(\varepsilon, \mathcal{F}, \mathcal{G}) \asymp \varepsilon^{q(\mathcal{F}, \mathcal{G})}.$$

Also set $q(\mathcal{G}) = q(\mathcal{G}, \mathcal{G})$.

Without loss of generality, assume $q(\mathcal{F}_1) \geq q(\mathcal{F}_2)$. Throughout the remainder of the paper $C$ is used to denote a generic constant which may vary from place to place and set $\mathcal{F} = \mathcal{F}_1 \cup \mathcal{F}_2$. Note that $q(\mathcal{F}_1, \mathcal{F}) = \min\{q(\mathcal{F}_1), q(\mathcal{F}_1, \mathcal{F}_2)\}$ and $q(\mathcal{F}) = \min\{q(\mathcal{F}_1), q(\mathcal{F}_2), q(\mathcal{F}_1, \mathcal{F}_2)\}$. In this setup strongly adaptive confidence intervals exist if and only if $q(\mathcal{F}_1, \mathcal{F}) = q(\mathcal{F}_1)$ or equivalently $q(\mathcal{F}_1) \leq q(\mathcal{F}_1, \mathcal{F}_2)$.

There are four cases of interest.

*Case* 1. $q(\mathcal{F}_2) \leq q(\mathcal{F}_1) \leq q(\mathcal{F}_1, \mathcal{F}_2)$. In this case $q(\mathcal{F}_1, \mathcal{F}) = q(\mathcal{F}_1)$ and strongly adaptive confidence intervals exist. These intervals have maximum expected length which can attain the same optimal rate of convergence as the minimax confidence interval over known $\mathcal{F}_i$. Specific shape restricted examples are given in Section 3.2 which illustrate this case and more general theory.

*Case* 2. $q(\mathcal{F}_1, \mathcal{F}_2) = q(\mathcal{F}_2) < q(\mathcal{F}_1)$. In this case $q(\mathcal{F}_1, \mathcal{F}) < q(\mathcal{F}_1)$ and thus strongly adaptive confidence intervals do not exist. Adaptive confidence intervals of level $1 - \alpha$ over $\mathcal{F}_1$ and $\mathcal{F}_2$ have maximum expected length over $\mathcal{F}_1$ which satisfies

$$(21) \qquad L(CI, \mathcal{F}_1) \geq \left(\frac{1}{2} - \alpha\right)\omega_+\left(\frac{z_\alpha}{\sqrt{n}}, \mathcal{F}_1, \mathcal{F}\right) \asymp n^{-q(\mathcal{F})/2}.$$

In contrast, if it is known that $f \in \mathcal{F}_1$, $1 - \alpha$ level confidence intervals can be constructed which satisfy

$$L(CI, \mathcal{F}_1) \leq Cn^{-q(\mathcal{F}_1)/2} \ll Cn^{-q(\mathcal{F})/2}.$$

Hence from this point of view the cost of adaptation is substantial. The rate of convergence of the maximum expected length of $CI$ over $\mathcal{F}_1$ is the same as that for the maximum expected length over $\mathcal{F}$.

EXAMPLE 1. Consider estimating the linear functional $Tf = f(0)$ over Lipschitz classes based on the Gaussian observations given in (1). For $0 < \beta \leq 1$ and $-\frac{1}{2} \leq a < b \leq \frac{1}{2}$, the Lipschitz function class over the interval $[a, b]$ is defined as

$$F(\beta, M, [a, b]) \tag{22}$$
$$= \{f : [-\tfrac{1}{2}, \tfrac{1}{2}] \to \mathbb{R}, |f(x) - f(y)| \leq M|x - y|^\beta \text{ for } x, y \in [a, b]\}.$$

It is also convenient to write $F(\beta, M)$ for $F(\beta, M, [-\tfrac{1}{2}, \tfrac{1}{2}])$.

Let $0 < \beta_2 < \beta_1 \leq 1$, set $\mathcal{F}_i = F(\beta_i, M)$ for $i = 1, 2$. In this case standard calculations as, for example, outlined in Cai and Low (2002) show that $\omega(\varepsilon, \mathcal{F}_1) = C\varepsilon^{2\beta_1/(2\beta_1+1)}(1 + o(1))$ and $\omega(\varepsilon, \mathcal{F}_1, \mathcal{F}_2) = C\varepsilon^{2\beta_2/(2\beta_2+1)}(1 + o(1))$. Hence

$$q(\mathcal{F}_1, \mathcal{F}) = q(\mathcal{F}_1, \mathcal{F}_2) = \frac{2\beta_2}{2\beta_2 + 1} < q(\mathcal{F}_1) = \frac{2\beta_1}{2\beta_1 + 1}.$$

*Case* 3. $q(\mathcal{F}_2) < q(\mathcal{F}_1, \mathcal{F}_2) < q(\mathcal{F}_1)$. In this case $q(\mathcal{F}_1, \mathcal{F}) < q(\mathcal{F}_1)$ and strongly adaptive confidence intervals do not exist. Any $1 - \alpha$ level adaptive confidence interval $CI$ over $\mathcal{F}_1$ and $\mathcal{F}_2$, must have maximum expected length of $CI$ over $\mathcal{F}_1$ satisfying

$$(23) \quad L(CI, \mathcal{F}_1) \geq \left(\frac{1}{2} - \alpha\right) \omega_+\left(\frac{z_\alpha}{\sqrt{n}}, \mathcal{F}_1, \mathcal{F}\right) \asymp n^{-q(\mathcal{F}_1, \mathcal{F}_2)/2} \gg n^{-q(\mathcal{F}_1)/2}.$$

The cost of adaptation in this case is that the rate of convergence of the maximum expected length of $CI$ over $\mathcal{F}_1$ is slower than that if it is known that $f \in \mathcal{F}_1$ but faster than for the maximum expected length over $\mathcal{F}_2$. An example for this case can be given as follows.

EXAMPLE 2. Suppose that the white noise with drift process (1) is observed and that the linear functional $Tf = f(0)$. Let the Lipschitz class $F(\beta, M, [a, b])$ be defined as above and let $\mathcal{D}$ be the set of all decreasing functions on $[-\tfrac{1}{2}, \tfrac{1}{2}]$. Set

$$F_D(\beta_1, M_1, \beta_2, M_2) = F(\beta_1, M_1, [-\tfrac{1}{2}, 0]) \cap F(\beta_2, M_2, [0, \tfrac{1}{2}]) \cap \mathcal{D}.$$

Let $\mathcal{F}_1 = F_D(\gamma_1, M_1, \gamma_2, M_2)$ and $\mathcal{F}_2 = F_D(\beta_1, N_1, \beta_2, N_2)$ with $1 \geq \gamma_1 > \gamma_2 > \beta_1 > \beta_2 > 0$. Then as in Cai and Low (2002) it is easy to check that

$$\omega(\varepsilon, \mathcal{F}_1) = C\varepsilon^{2\gamma_1/(2\gamma_1+1)}(1 + o(1)),$$
$$\omega(\varepsilon, \mathcal{F}_2) = C\varepsilon^{2\beta_1/(2\beta_1+1)}(1 + o(1)),$$
$$(24)$$
$$\omega(\varepsilon, \mathcal{F}_1, \mathcal{F}_2) = C\varepsilon^{2\gamma_2/(2\gamma_2+1)}(1 + o(1)),$$
$$\omega(\varepsilon, \mathcal{F}_2, \mathcal{F}_1) = C\varepsilon^{2\gamma_1/(2\gamma_1+1)}(1 + o(1)).$$

Note that in this case $\omega(\varepsilon, \mathcal{F}_1, \mathcal{F}_2) \neq \omega(\varepsilon, \mathcal{F}_2, \mathcal{F}_1)(1+o(1))$. Since $\gamma_1 > \gamma_2$, it then follows from (24) that

$$q(\mathcal{F}_1, \mathcal{F}) = q(\mathcal{F}_1, \mathcal{F}_2) = \frac{2\gamma_2}{2\gamma_2 + 1}.$$

Hence $0 < q(\mathcal{F}_2) < q(\mathcal{F}_1, \mathcal{F}) < q(\mathcal{F}_1) < 1$.

*Case* 4. $q(\mathcal{F}_1, \mathcal{F}_2) < q(\mathcal{F}_2) \leq q(\mathcal{F}_1)$. In this case, strongly adaptive confidence intervals do not exist and the cost of adaptation is extraordinary. If $f$ is known to be in $\mathcal{F}_i$, one can attain the rate of convergence $n^{q(\mathcal{F}_i)/2}$ for the maximum expected length of the optimal $1 - \alpha$ level confidence interval over $\mathcal{F}_i$. Without the information $1 - \alpha$ level adaptive confidence intervals over $\mathcal{F}_1$ and $\mathcal{F}_2$ must have maximum expected length over $\mathcal{F}_i$ at least of order $n^{-q(\mathcal{F}_1, \mathcal{F}_2)/2}$. An example is given below.

EXAMPLE 3. Once again consider the white noise model with $Tf = f(0)$. Let

$$F(\beta_1, M_1, \beta_2, M_2) = F(\beta_1, M_1, [-\tfrac{1}{2}, 0]) \cap F(\beta_2, M_2, [0, \tfrac{1}{2}])$$

and consider $0 < \gamma_2 \leq \gamma_1 \leq 1$ and $0 < \beta_1 \leq \beta_2 \leq 1$. Set $\mathcal{F}_1 = F(\gamma_1, M_1, \gamma_2, M_2)$ and $\mathcal{F}_2 = F(\beta_1, N_1, \beta_2, N_2)$. Standard calculations show that $\omega(\varepsilon, \mathcal{F}_1) = C\varepsilon^{2\gamma_1/(2\gamma_1+1)}(1+o(1))$ and $\omega(\varepsilon, \mathcal{F}_2) = C\varepsilon^{2\beta_2/(2\beta_2+1)}(1+o(1))$. The between class modulus is given as

$$(25) \qquad \omega(\varepsilon, \mathcal{F}_1, \mathcal{F}_2) = C\varepsilon^{2\rho/(2\rho+1)}(1+o(1))$$

where $\rho = \max(\min(\gamma_1, \beta_1), \min(\gamma_2, \beta_2))$.

When $\gamma_1 \geq \beta_2 > \beta_1 \geq \gamma_2$, the quantity $\rho$ in (25) equals $\beta_1$ and hence

$$q(\mathcal{F}_1, \mathcal{F}_2) = \frac{2\beta_1}{2\beta_1 + 1}.$$

Therefore in this case $q(\mathcal{F}_1, \mathcal{F}_2) < \min(q(\mathcal{F}_1), q(\mathcal{F}_2))$.

**3. Adaptation over nested parameter spaces.** Section 2 gave the adaptation theory for two convex parameter spaces. This theory can be extended to more general collections of parameter spaces. In this section the focus is on adaptation over a collection of a finite or countably infinite number of nested convex parameter spaces, $\mathcal{F}_1 \subset \mathcal{F}_2 \subset \cdots \subset \mathcal{F}_k$, where in the case of $k = \infty$, $\mathcal{F}_\infty$ denotes $\bigcup_{i=1}^\infty \mathcal{F}_i$. The objective is, for a given linear functional $Tf$, to construct variable length confidence intervals which have coverage probability of at least $1 - \alpha$ over $\mathcal{F}_k$ and which simultaneously minimize the expected length over each of the parameter spaces $\mathcal{F}_j$. A target for these expected lengths has been provided by the lower bound given in Theorem 1, namely

$$
(26) \quad L^*_\alpha(\mathcal{F}_j, \mathcal{F}_k) \geq \left(\frac{1}{2} - \alpha\right) \omega_+\left(\frac{z_\alpha}{\sqrt{n}}, \mathcal{F}_j, \mathcal{F}_k\right)
$$

where $\omega_+(\varepsilon, \mathcal{F}_j, \mathcal{F}_k)$ is the between class modulus as given in (5).

The major result of this section is to show that adaptive confidence intervals exist and to construct such adaptive intervals. As in Section 2.2 the construction of these adaptive confidence procedures relies on the ordered modulus $\omega(\varepsilon, \mathcal{F}_i, \mathcal{F}_j)$ as given in (4). For $1 \leq i, j \leq k$ set $\omega_{i,j} = \omega(\frac{z_{\alpha/2}}{\sqrt{n}}, \mathcal{F}_i, \mathcal{F}_j)$ and let $\hat{T}_{i,j}$ be linear estimators with variances and biases bounded as in (13)–(15).

The confidence procedure is built in two steps. In the first step for each $1 \leq j \leq k$ an interval is constructed which controls the coverage probability over $\mathcal{F}_k$ and which also has expected length over $\mathcal{F}_j$ within a constant factor of the lower bound given by (26). In the second step these intervals are combined to create a single interval which maintains coverage while simultaneously attaining an expected length over every $\mathcal{F}_j$ within a fixed constant factor of the lower bound given in (26).

For the first step define the confidence intervals $CI^*_j$ as follows. For $1 \leq j \leq k$ set $\xi_j = \omega_+(\frac{z_{\alpha/2}}{\sqrt{n}}, \mathcal{F}_j, \mathcal{F}_k)$ and define $CI^*_j$ by

$$
(27) \quad \begin{aligned} CI^*_j = \Big[ & \frac{\hat{T}_{j,k} + \hat{T}_{k,j}}{2} - \{(\hat{T}_{j,k} - \hat{T}_{k,j})_+ + 2\xi_j\}, \\ & \frac{\hat{T}_{j,k} + \hat{T}_{k,j}}{2} + \{(\hat{T}_{j,k} - \hat{T}_{k,j})_+ + 2\xi_j\} \Big]. \end{aligned}
$$

Lemma 2 shows that these intervals have guaranteed coverage over $\mathcal{F}_k$ and near optimal expected length over $\mathcal{F}_j$.

REMARK 4. This interval is designed for $0 < \alpha \leq 0.2$. If $0.2 < \alpha \leq 0.5$ all subsequent results hold with minor modifications, as noted in later remarks, when the interval is replaced by

$$
(28) \quad \begin{aligned} CI^*_j = \Big[ & \frac{\hat{T}_{j,k} + \hat{T}_{k,j}}{2} - \{(\hat{T}_{j,k} - \hat{T}_{k,j})_+ + 3\xi_j\}, \\ & \frac{\hat{T}_{j,k} + \hat{T}_{k,j}}{2} + \{(\hat{T}_{j,k} - \hat{T}_{k,j})_+ + 3\xi_j\} \Big]. \end{aligned}
$$

LEMMA 2. *For $0 < \alpha \leq 0.2$, the confidence interval $CI_j^*$ defined in (27) has coverage probability of at least $1 - \frac{2}{7}\alpha$ for all $f \in \mathcal{F}_k$ and satisfies*

$$
\begin{aligned}
L(CI_j^*, \mathcal{F}_j) &\leq \left\{ 2\Phi\left(\frac{1}{2}z_{\alpha/2}\right) + \frac{4}{\sqrt{2\pi}z_{\alpha/2}} \exp\left(-\frac{1}{8}z_{\alpha/2}^2\right) + 4 \right\} \cdot \xi_j \\
&\leq 8\omega_+\left(\frac{z_{\alpha/2}}{\sqrt{n}}, \mathcal{F}_j, \mathcal{F}_k\right).
\end{aligned}
$$
(29)

REMARK 5. For $0.2 < \alpha \leq 0.5$ the interval given in (28) satisfies the same coverage but has expected length bounded by $10\omega_+(\frac{z_{\alpha/2}}{\sqrt{n}}, \mathcal{F}_j, \mathcal{F}_k)$.

In the following proof, and throughout the rest of the paper, write $Z$ for a standard Normal random variable.

PROOF OF LEMMA 2. Lemma 2 gives a bound on both coverage probability and expected length. First consider coverage probability. It is easy to see that the interval $CI_j^*$ contains the interval $CI_j$ defined as

$$CI_j = [\hat{T}_{k,j} - 2\xi_j, \hat{T}_{j,k} + 2\xi_j] \tag{30}$$

where the interval $CI_j$ is taken to be the empty set whenever the left endpoint of the above interval is larger than the right endpoint. First note that for $f \in \mathcal{F}_k$, $E\hat{T}_{k,j} - Tf \leq \frac{1}{2}\omega_{k,j}$ and that $E\hat{T}_{j,k} - Tf \geq -\frac{1}{2}\omega_{j,k}$. Let

$$z_{k,j} = \frac{\hat{T}_{k,j} - Tf - (1/2)\omega_{k,j}}{\omega_{k,j}/z_{\alpha/2}},$$

$$z_{j,k} = \frac{\hat{T}_{j,k} - Tf + (1/2)\omega_{j,k}}{\omega_{j,k}/z_{\alpha/2}}.$$

Then for any $f \in \mathcal{F}_k$ it follows from (14) and (15) that $z_{k,j}$ has a Normal distribution with mean less than or equal to zero and variance bounded by 1, and $z_{j,k}$ has a Normal distribution with mean greater than or equal to zero and variance bounded by 1. Note that $\xi_j = \max(\omega_{k,j}, \omega_{j,k})$. Hence for $f \in \mathcal{F}_k$,

$$
\begin{aligned}
P(Tf \notin CI_j^*) &\leq P(Tf \notin CI_j) \\
&\leq P\left(z_{k,j} > \left(2\frac{\xi_j}{\omega_{k,j}} - \frac{1}{2}\right)z_{\alpha/2}\right) + P\left(z_{j,k} \geq \left(-2\frac{\xi_j}{\omega_{j,k}} + \frac{1}{2}\right)z_{\alpha/2}\right) \\
&\leq 2P\left(Z \geq \frac{3}{2}z_{\alpha/2}\right).
\end{aligned}
$$

Note that for a fixed $\lambda > 1$, it is easy to verify that $g(z) = P(Z \geq \lambda z)/P(Z \geq z)$ is a strictly decreasing function of $z$ for $z > 0$ and for $\alpha = 0.2$,

$$2P(Z \geq \tfrac{3}{2}z_{\alpha/2}) \leq \tfrac{2}{7}\alpha.$$

Hence, $P(Tf \notin CI_j^*) \leq \frac{2}{7}\alpha$ and so the claim of the required coverage probability has been established.

Now turn to the bound on expected length given in (29) for which the following technical lemma is needed.

LEMMA 3. *Let $X \sim N(\mu, \sigma^2)$ with $\mu \leq \mu_0$ and $0 < \sigma \leq \sigma_0$. Then*

$$(31) \qquad EX\mathbb{1}(X > 0) \leq \mu_0 \Phi\left(\frac{\mu_0}{\sigma_0}\right) + \frac{\sigma_0}{\sqrt{2\pi}} \exp\left(-\frac{\mu_0^2}{2\sigma_0^2}\right).$$

PROOF. It is easy to check by taking partial derivatives that $EX\mathbb{1}(X > 0)$ is an increasing function of both $\mu$ and $\sigma$. Hence

$$E_{\mu,\sigma} X\mathbb{1}(X > 0)$$
$$\leq E_{\mu_0,\sigma_0} X\mathbb{1}(X > 0)$$
$$= \frac{1}{\sqrt{2\pi}\sigma_0} \int_0^\infty x \exp\left(-\frac{(x-\mu_0)^2}{2\sigma_0^2}\right) dx$$
$$= \frac{1}{\sqrt{2\pi}\sigma_0} \int_0^\infty \mu_0 \exp\left(-\frac{(x-\mu_0)^2}{2\sigma_0^2}\right) dx + \frac{\sigma_0}{\sqrt{2\pi}} \int_{-\mu_0/\sigma_0}^\infty y \exp\left(-\frac{y^2}{2}\right) dy$$
$$= \mu_0 \Phi\left(\frac{\mu_0}{\sigma_0}\right) + \frac{\sigma_0}{\sqrt{2\pi}} \exp\left(-\frac{\mu_0^2}{2\sigma_0^2}\right). \qquad \square$$

Now note that for $f \in \mathcal{F}_j$,

$$E(\hat{T}_{j,k} - \hat{T}_{k,j}) \leq \xi_j \quad \text{and} \quad \text{Var}(\hat{T}_{j,k} - \hat{T}_{k,j}) \leq \frac{4}{z_{\alpha/2}^2}\xi_j^2,$$

and so from Lemma 3 it follows that

$$(32) \quad E(\hat{T}_{j,k} - \hat{T}_{k,j})_+ \leq \left\{\Phi\left(\frac{1}{2}z_{\alpha/2}\right) + \frac{2}{\sqrt{2\pi}z_{\alpha/2}} \exp\left(-\frac{1}{8}z_{\alpha/2}^2\right)\right\} \cdot \xi_j \leq 2\xi_j$$

and hence (29) is satisfied. $\square$

Lemma 2 shows that the interval $CI_j^*$ has guaranteed coverage over $\mathcal{F}_k$ and near optimal expected length over $\mathcal{F}_j$. Before turning to the construction of an adaptive confidence interval we state a simple preliminary lemma. The proof is straightforward and not given here.

LEMMA 4. *Let $0 < \xi_1 \leq \xi_2 \leq \cdots \leq \xi_k$ be a sequence of monotonically increasing positive numbers. Then there exists a unique subsequence $\xi_{j_1} < \xi_{j_2} < \cdots < \xi_{j_m}$ with $j_m = k$, such that for all $1 \leq i \leq m$,*

$$(33) \qquad \xi_{j_i} \geq 2\xi_{j_{i-1}} \quad \text{and} \quad \xi_{j_i} < 2\xi_j \qquad \text{for all } j_{i-1} < j < j_i$$

*where we set $j_0 = 0$ and $\xi_0 = 0$.*

The construction of the adaptive confidence interval proceeds as follows. Once again for $1 \leq j \leq k$, set $\xi_j = \omega_+(\frac{z_{\alpha/2}}{\sqrt{n}}, \mathcal{F}_j, \mathcal{F}_k)$. Let $\xi_{j_1} < \xi_{j_2} < \cdots < \xi_{j_m}$ be the subsequence satisfying (33). Let $\hat{j}$ be the index of the shortest interval among all the $CI^*_{j_i}$ for $1 \leq i \leq m$. More precisely,

$$\hat{j} = \operatorname*{arg\,min}_{j_i, 1 \leq i \leq m} L(CI^*_{j_i}).$$

Then the adaptive confidence interval for $Tf$ is defined by

(34) $$CI^* = CI^*_{\hat{j}}.$$

The following theorem shows that $CI^*$ is a $1-\alpha$ level adaptive confidence interval over the collection $\{\mathcal{F}_j, j = 1, \ldots, k\}$.

THEOREM 2. *The confidence interval $CI^*$ defined in (34) has coverage probability of at least $1 - \alpha$ for all $f \in \mathcal{F}_k$, that is, $CI^* \in \mathcal{I}_{\alpha, \mathcal{F}_k}$ and satisfies*

(35) $$L^*_\alpha(\mathcal{F}_j, \mathcal{F}_k) \leq L(CI^*, \mathcal{F}_j) \leq \frac{16 z_{\alpha/2}}{(1/2 - \alpha) z_\alpha} L^*_\alpha(\mathcal{F}_j, \mathcal{F}_k)$$

*simultaneously for all $1 \leq j \leq k$. Moreover,*

(36) $$L(CI^*, \mathcal{F}_{j_i}) \leq 8 \omega_+\left(\frac{z_{\alpha/2}}{\sqrt{n}}, \mathcal{F}_{j_i}, \mathcal{F}_k\right)$$

*for all $1 \leq i \leq m$, and for all $1 \leq j \leq k$*

(37) $$L(CI^*, \mathcal{F}_j) \leq 16 \omega_+\left(\frac{z_{\alpha/2}}{\sqrt{n}}, \mathcal{F}_j, \mathcal{F}_k\right).$$

The proof of Theorem 2 rests on the following important technical lemma. Recall that Lemma 2 gives a lower bound on coverage over $\mathcal{F}_k$ and an upper bound on expected length over $\mathcal{F}_j$. Lemma 5 shows, in a precise way, that if $CI^*_j$ has a large expected length it must have high coverage probability.

LEMMA 5. *If $f \in \mathcal{F}_k$ and*

$$P(Tf \notin CI^*_j) > 2P(Z \geq \tfrac{1}{4}(\lambda + 3) z_{\alpha/2}),$$

*then*

$$E(\hat{T}_{j,k} - \hat{T}_{k,j}) \leq \lambda \xi_j.$$

PROOF. First note that

$$P(Tf \notin CI^*_j) \leq P\left(Tf \leq \frac{\hat{T}_{j,k} + \hat{T}_{k,j}}{2} - (\hat{T}_{j,k} - \hat{T}_{k,j} + 2\xi_j)\right)$$
$$+ P\left(Tf \geq \frac{\hat{T}_{j,k} + \hat{T}_{k,j}}{2} + (\hat{T}_{j,k} - \hat{T}_{k,j} + 2\xi_j)\right).$$

Now note that

$$-\frac{1}{2}\xi_j - \frac{1}{2}E(\hat{T}_{j,k} - \hat{T}_{k,j}) \leq E\frac{\hat{T}_{j,k} + \hat{T}_{k,j}}{2} - Tf \leq \frac{1}{2}\xi_j + \frac{1}{2}E(\hat{T}_{j,k} - \hat{T}_{k,j}).$$

Let $X = \frac{\hat{T}_{j,k} + \hat{T}_{k,j}}{2} - Tf - (\hat{T}_{j,k} - \hat{T}_{k,j} + 2\xi_j)$. Suppose that

$$E(\hat{T}_{j,k} - \hat{T}_{k,j}) > \lambda \xi_j.$$

Then

$$E(X) \leq -\tfrac{1}{2}(\lambda + 3)\xi_j \quad \text{and} \quad \text{Var}(X) \leq \frac{4}{z_{\alpha/2}^2}\xi_j^2.$$

Hence

$$P\left(Tf \leq \frac{\hat{T}_{j,k} + \hat{T}_{k,j}}{2} - (\hat{T}_{j,k} - \hat{T}_{k,j} + 2\xi_j)\right) = P(X \geq 0)$$
$$\leq P\left(Z \geq \frac{1}{4}(\lambda + 3)z_{\alpha/2}\right).$$

Similarly,

$$P\left(Tf \geq \frac{\hat{T}_{j,k} + \hat{T}_{k,j}}{2} + (\hat{T}_{j,k} - \hat{T}_{k,j} + 2\xi_j)\right) \leq P\left(Z \geq \frac{1}{4}(\lambda + 3)z_{\alpha/2}\right).$$

Hence,

$$P(Tf \notin CI_j^*) \leq 2P(Z \geq \tfrac{1}{4}(\lambda + 3)z_{\alpha/2}). \qquad \square$$

PROOF OF THEOREM 2. Note that it suffices to prove (36) since (37) follows immediately from (36) and (35) is a direct consequence of (20), (7) and (37). For (36) assume without loss of generality that $\xi_j \geq 2\xi_{j-1}$ for all $1 \leq j \leq k$; otherwise we can work along the subsequence. First note that since $CI^*$ is the shortest of all the $CI_j^*$ confidence intervals Lemma 2 yields that the expected length of $CI^*$ satisfies

$$L(CI^*, \mathcal{F}_j) \leq L(CI_j^*, \mathcal{F}_j)$$
(38)
$$\leq \left\{2\Phi\left(\frac{1}{2}z_{\alpha/2}\right) + \frac{4}{\sqrt{2\pi}z_{\alpha/2}}\exp\left(-\frac{1}{8}z_{\alpha/2}^2\right) + 4\right\} \cdot \xi_j$$
$$\leq 8\xi_j.$$

Now turn to the proof of coverage. Note that

(39)
$$P(Tf \notin CI^*) = \sum_{j=1}^{k} P(Tf \notin CI_j^* \cap \hat{j} = j)$$
$$\leq \sum_{j=1}^{k} \min\{P(Tf \notin CI_j^*), P(\hat{j} = j)\}.$$

For $l \geq 0$, denote $d(l) = 2P(Z \geq \frac{1}{4}(l+6)z_{\alpha/2})$. Note that $d(0) = 2P(Z \geq \frac{3}{2}z_{\alpha/2}) \leq \frac{2}{7}\alpha$. For $l \geq 1$ let

$$(40) \qquad A_l = \{j : d(l) < P(Tf \notin CI_j^*) \leq d(l-1)\}$$

and let $j(l) = \min\{j : j \in A_l\}$. Note that it follows from Lemma 2 that $\bigcup_l A_l = \{j \geq 1\}$. Then by Lemma 5

$$(41) \qquad E(\hat{T}_{j(l),k} - \hat{T}_{k,j(l)}) \leq (l+3)\xi_{j(l)}.$$

Note that $\text{Var}(\hat{T}_{j(l),k} - \hat{T}_{k,j(l)}) \leq \frac{4}{z_{\alpha/2}^2}\xi_j^2$, so

$$P(L(CI_{j(l)}^*) > 4\rho\xi_{j(l)}) = P(\hat{T}_{j(l),k} - \hat{T}_{k,j(l)} > 2(\rho-1)\xi_{j(l)})$$
$$\leq P(Z \geq (\rho - \tfrac{5}{2} - \tfrac{1}{2}l)z_{\alpha/2}).$$

Since $\xi_j \geq 2\xi_{j-1}$, it follows that, for any integer $m > 0$,

$$P(\hat{j} \geq j(l) + m) \leq P(L(CI_{j(l)}^*) > 4\xi_{j(l)+m})$$
$$\leq P(L(CI_{j(l)}^*) > 4 \cdot 2^m \xi_{j(l)})$$
$$\leq P(Z \geq (2^m - \tfrac{5}{2} - \tfrac{1}{2}l)z_{\alpha/2})$$
$$\equiv \gamma_{l,m}.$$

Let $j_* = \min\{j(l) : 1 \leq l \leq 8\}$. For $m = 3$ and $1 \leq l \leq 8$, $\gamma_{l,m} = P(Z \geq \frac{1}{2}(11-l)z_{\alpha/2})$. If $j_* = j(1)$, then

$$\sum_{j=j_*}^{k} \min\{P(Tf \notin CI_j^*), P(\hat{j}=j)\} \leq d(0) + d(0) + d(0) + \gamma_{1,3}$$
$$\leq \tfrac{6}{7}\alpha + P(Z \geq 5z_{\alpha/2}).$$

Similarly, if $j_* = j(l)$ for some $2 \leq l \leq 8$, then

$$\sum_{j=j_*}^{k} \min\{P(Tf \notin CI_j^*), P(\hat{j}=j)\} \leq d(l-1) + d(0) + d(0) + \gamma_{l,3}$$
$$\leq \tfrac{6}{7}\alpha + P(Z \geq 5z_{\alpha/2}).$$

Hence

$$(42) \qquad \sum_{j=j_*}^{k} \min\{P(Tf \notin CI_j^*), P(\hat{j}=j)\} \leq \tfrac{6}{7}\alpha + P(Z \geq 5z_{\alpha/2}).$$

The following simple lemma can be used to bound $P(Z \geq 5z_{\alpha/2})$.

LEMMA 6. *Let $Z \sim N(0,1)$ and let $a > 0$ and $b > 0$ be two constants. Then*

$$P(Z \geq a+b) \leq \exp(-(ab + \tfrac{1}{2}b^2))P(Z > a).$$

Applying Lemma 6 with $a = z_{\alpha/2}$ and $b = 4z_{\alpha/2}$, it follows that

$$P(Z \geq 5z_{\alpha/2}) = P(Z \geq z_{\alpha/2} + 4z_{\alpha/2})$$
$$\leq \exp(-12z_{\alpha/2}^2) \cdot \frac{\alpha}{2} \leq \frac{1}{14}\alpha.$$

Therefore

$$P(Tf \notin CI^*) \leq \sum_{j=1}^{k} \min\{P(Tf \notin CI_j^*), P(\hat{j} = j)\}$$
$$\leq \tfrac{13}{14}\alpha + \sum_{l=9}^{\infty} \sum_{j \in A_l} \min\{P(Tf \notin CI_j^*), P(\hat{j} = j)\}.$$

For $l \geq 9$, let $m_l$ be the smallest integer satisfying $2^{m_l} \geq \tfrac{1}{4}(3l+7)$. Then $m_l \leq \log_2(3l+7) - 1$. Recall that for $j \in A_l$, $P(Tf \notin CI_j^*) \leq 2P(Z \geq \tfrac{1}{4}(l+3)z_{\alpha/2})$. Now note that

$$P(\hat{j} \geq j(l) + m_l) \leq \gamma_{l,m_l} \leq P(Z \geq \tfrac{1}{4}(l+3)z_{\alpha/2}).$$

So, for $l \geq 9$,

$$\sum_{j \in A_l} \min\{P(Tf \notin CI_j^*), P(\hat{j} = j)\}$$
$$\leq m_l \cdot 2P(Z \geq \tfrac{1}{4}(l+3)z_{\alpha/2}) + \gamma_{l,m_l}$$
$$\leq (2m_l + 1)P(Z \geq \tfrac{1}{4}(l+3)z_{\alpha/2}).$$

So

$$\sum_{l=9}^{\infty} \sum_{j \in A_l} \min\{P(Tf \notin CI_j^*), P(\hat{j} = j)\}$$
$$\leq \sum_{l=9}^{\infty} (2\log_2(3l+7) - 1)P(Z \geq \tfrac{1}{4}(l+3)z_{\alpha/2}).$$

Lemma 6 yields

$$P\left(Z \geq \frac{1}{4}(l+3)z_{\alpha/2}\right) \leq P\left(Z \geq z_{\alpha/2} + \frac{1}{4}(l-1)z_{\alpha/2}\right)$$
$$\leq \exp\left(-\left(\frac{1}{4}(l-1) + \frac{1}{32}(l-1)^2\right)z_{\alpha/2}^2\right) \cdot \frac{\alpha}{2}.$$

Hence,

$$\sum_{l=9}^{\infty} \sum_{j \in A_l} \min\{P(Tf \notin CI_j^*), P(\hat{j} = j)\}$$

$$\leq \frac{\alpha}{2} \sum_{l=9}^{\infty} (2\log_2(3l+7) - 1) \exp\left(-\left(\frac{1}{4}(l-1) + \frac{1}{32}(l-1)^2\right) z_{\alpha/2}^2\right)$$

$$\leq \frac{\alpha}{2} \sum_{l=8}^{\infty} (2\log_2(3l+10) - 1) \exp\left(-\frac{z_{\alpha/2}^2}{32} l^2\right) \exp\left(-\frac{z_{\alpha/2}^2}{4} l\right).$$

It is easy to see that for $l \geq 8$, $(2\log_2(3l+10) - 1) \exp(-(z_{\alpha/2}^2/32)l^2)$ is strictly decreasing and

$$(2\log_2(3l+10) - 1) \exp\left(-\frac{z_{\alpha/2}^2}{32} l^2\right) \leq \frac{1}{2}.$$

So,

$$\sum_{l=9}^{\infty} \sum_{j \in A_l} \min\{P(Tf \notin CI_j^*), P(\hat{j} = j)\} \leq \frac{1}{4}\alpha \sum_{l=8}^{\infty} \exp\left(-\frac{z_{\alpha/2}^2}{4} l\right) \leq \frac{1}{14}\alpha.$$

Hence,

$$P(Tf \notin CI^*) \leq \tfrac{13}{14}\alpha + \tfrac{1}{14}\alpha = \alpha. \qquad \Box$$

3.1. *Adaptation over nearly nested parameter spaces.* In some common cases of interest such as Hölder spaces, Sobolev spaces and Besov spaces, the parameter spaces are not exactly nested, but have nested structure in terms of the moduli of continuity. Theorem 2 can be generalized to such nearly nested parameter spaces.

Denote by C.Hull($\mathcal{F}$) the convex hull of a parameter set $\mathcal{F}$. Let $\mathcal{F}_i$, $i = 1, \ldots, k$, be convex parameter spaces and for any integer $1 \leq m \leq k$ let $\mathcal{G}_m = \bigcup_{l=1}^{m} \mathcal{F}_l$. Suppose the following condition, which is trivially satisfied if $\mathcal{F}_i$ are nested, holds.

CONDITION C. For $1 \leq j \leq k$ and some constants $C_2 \geq C_1 > 0$,

$$\omega(\varepsilon, \text{C.Hull}(\mathcal{G}_j), \text{C.Hull}(\mathcal{G}_k)) \leq C_1 \omega(\varepsilon, \mathcal{G}_j, \mathcal{G}_k) \leq C_2 \omega(\varepsilon, \mathcal{F}_j, \mathcal{F}_k)$$

and

$$\omega(\varepsilon, \text{C.Hull}(\mathcal{G}_k), \text{C.Hull}(\mathcal{G}_j)) \leq C_1 \omega(\varepsilon, \mathcal{G}_k, \mathcal{G}_j) \leq C_2 \omega(\varepsilon, \mathcal{F}_k, \mathcal{F}_j)$$

for all $0 < \varepsilon < \varepsilon_0$.

Similarly to the nested case for $1 \leq i, j \leq k$, set $\omega'_{i,j} = \omega(\frac{z_{\alpha/2}}{\sqrt{n}}, \text{C.Hull}(\mathcal{G}_i), \text{C.Hull}(\mathcal{G}_j))$, and once again Cai and Low (2004) give a construction of linear estimators $\hat{T}'_{i,j}$ which have variance bounded by

$$\text{Var}(\hat{T}'_{i,j}) \leq \frac{1}{z^2_{\alpha/2}} {\omega'_{i,j}}^2$$

and bias which satisfies

$$\inf_{f \in \mathcal{F}_j} (E(\hat{T}'_{i,j}) - Tf) \geq -\tfrac{1}{2}\omega'_{i,j}$$

and

$$\sup_{f \in \mathcal{F}_i} (E(\hat{T}'_{i,j}) - Tf) \leq \tfrac{1}{2}\omega'_{i,j}.$$

Set $\xi'_j = \omega_+(\frac{z_{\alpha/2}}{\sqrt{n}}, \text{C.Hull}(\mathcal{G}_j), \text{C.Hull}(\mathcal{G}_k))$ and define the confidence intervals $CI^*_j$ as earlier. When $0 < \alpha \leq 0.2$, let

(43)
$$CI^*_j = \Big[\frac{\hat{T}'_{j,k} + \hat{T}'_{k,j}}{2} - \{(\hat{T}'_{j,k} - \hat{T}'_{k,j})_+ + 2\xi'_j\},$$
$$\frac{\hat{T}'_{j,k} + \hat{T}'_{k,j}}{2} + \{(\hat{T}'_{j,k} - \hat{T}'_{k,j})_+ + 2\xi'_j\}\Big]$$

and when $0.2 < \alpha < 0.5$ let

(44)
$$CI^*_j = \Big[\frac{\hat{T}'_{j,k} + \hat{T}'_{k,j}}{2} - \{(\hat{T}'_{j,k} - \hat{T}'_{k,j})_+ + 3\xi'_j\},$$
$$\frac{\hat{T}'_{j,k} + \hat{T}'_{k,j}}{2} + \{(\hat{T}'_{j,k} - \hat{T}'_{k,j})_+ + 3\xi'_j\}\Big].$$

Following the argument given in the nested case let $\xi'_{j_i}$ be a subsequence of $\xi'_j$ satisfying (33) and let $\hat{j} = \arg\min_{j_i, 1 \leq i \leq m} L(CI^*_{j_i})$ be the index of the shortest interval along the subsequence and define the adaptive confidence interval for $Tf$ by

(45) $$CI^* = CI^*_{\hat{j}}.$$

As stated precisely in the following result this confidence interval is adaptive over the parameter spaces $\{\mathcal{F}_j : j = 1, \ldots, k\}$.

PROPOSITION 3. *Suppose Condition* C *holds. Then the confidence interval* $CI^*$ *defined in* (45) *has coverage probability of at least* $1 - \alpha$ *for all* $f \in \mathcal{F} = \bigcup_{j=1}^k \mathcal{F}_j$ *and satisfies the lower bound on expected length,*

(46) $$L^*_\alpha(\mathcal{F}_j, \mathcal{F}) \leq L(CI^*, \mathcal{F}_j) \leq C(\alpha) L^*_\alpha(\mathcal{F}_j, \mathcal{F}),$$

*simultaneously for all* $1 \leq j \leq k$, *where the constant* $C(\alpha)$ *only depends on* $\alpha$ *and is independent of* $k$. *In other words,* $L(CI^*, \mathcal{F}_j) \asymp \omega_+(\frac{z_{\alpha/2}}{\sqrt{n}}, \mathcal{F}_j, \mathcal{F})$ *for all* $1 \leq j \leq k$.

We omit the proof of Proposition 3 since it essentially follows a similar path to that of Theorem 2.

3.2. *Examples.* Theorem 2 and Proposition 3 have established general adaptation results for collections of nested or nearly nested parameter spaces. In this section a couple of examples are given which illustrate this general theory.

Suppose that we observe the white noise with drift process (1) and that the linear functional is point evaluation. For convenience take $Tf = f(0)$. Let $\mathcal{D}$ be the set of all decreasing functions on $[-\frac{1}{2}, \frac{1}{2}]$ and let $F_D(\beta, M) = F(\beta, M) \cap \mathcal{D}$ be the collection of monotonically decreasing Lipschitz functions where $F(\beta, M)$ is the Lipschitz class defined in (22).

For integer $j \geq 1$ let $M_j = 2^{j(2\beta+1)} \frac{1}{\sqrt{n}}$ and let $\mathcal{G} = \bigcup_{j=1}^\infty F_D(\beta, M_j)$. Standard calculations as in, for example, Donoho and Liu (1987), yield

$$(47) \quad \begin{aligned} \omega(\varepsilon, F_D(\beta, M), \mathcal{G}) &= \omega(\varepsilon, \mathcal{G}, F_D(\beta, M)) \\ &= (2\beta + 1)^{1/(2\beta+1)} M^{1/(2\beta+1)} \varepsilon^{2\beta/(2\beta+1)} \end{aligned}$$

for $M \geq (2\beta + 1)^{1/2} \varepsilon$. Let $\xi_j = \omega(\frac{z_{\alpha/2}}{\sqrt{n}}, F_D(\beta, M_j), \mathcal{G})$. Then it is easy to see that $\xi_{j+1} = 2\xi_j$ and hence the adaptive confidence interval given in (34) has coverage probability over $\mathcal{G}$ of at least $1 - \alpha$ and satisfies

$$(48) \quad L(CI^*, F_D(\beta, M_j)) \leq 6(2\beta + 1)^{1/(2\beta+1)} M_j^{1/(2\beta+1)} z_{\alpha/2}^{2\beta/(2\beta+1)} n^{-\beta/(2\beta+1)}.$$

Furthermore, for any $M > 0$,

$$(49) \quad L(CI^*, F_D(\beta, M)) \leq 12(2\beta + 1)^{1/(2\beta+1)} M^{1/(2\beta+1)} z_{\alpha/2}^{2\beta/(2\beta+1)} n^{-\beta/(2\beta+1)}$$

for all sufficiently large $n$.

Another common problem in function estimation is to adapt over smoothness classes. For fixed $M > 0$, the classes $F_D(\gamma_1, M) \subset F_D(\gamma_2, M)$ whenever $0 < \gamma_2 < \gamma_1 \leq 1$. Let $\mathcal{G}' = \bigcup_{0 \leq \gamma \leq 1} F_D(\gamma, M)$. Then once again standard calculations yield

$$(50) \quad \begin{aligned} \omega(\varepsilon, F_D(\beta, M), \mathcal{G}') &= \omega(\varepsilon, \mathcal{G}', F_D(\beta, M)) \\ &= (2\beta + 1)^{1/(2\beta+1)} M^{1/(2\beta+1)} \varepsilon^{2\beta/(2\beta+1)}. \end{aligned}$$

Now let $1 = \beta_1 > \beta_2 > \cdots$ be the sequence such that

$$\omega_+\left(\frac{z_{\alpha/2}}{\sqrt{n}}, F_D(\beta_{j+1}, M), \mathcal{G}'\right) = 2\omega_+\left(\frac{z_{\alpha/2}}{\sqrt{n}}, F_D(\beta_j, M), \mathcal{G}'\right).$$

Then the adaptive confidence interval given in (34) has coverage probability over $\mathcal{G}'$ of at least $1 - \alpha$ and satisfies

$$(51) \quad \begin{aligned} &L(CI^*, F_D(\beta_j, M)) \\ &\qquad \leq 6(2\beta_j + 1)^{1/(2\beta_j+1)} M^{1/(2\beta_j+1)} z_{\alpha/2}^{2\beta_j/(2\beta_j+1)} n^{-\beta_j/(2\beta_j+1)}. \end{aligned}$$

Furthermore, for any $0 < \beta \leq 1$,

(52) $\quad L(CI^*, F_D(\beta, M)) \leq 12(2\beta+1)^{1/(2\beta+1)} M^{1/(2\beta+1)} z_{\alpha/2}^{2\beta/(2\beta+1)} n^{-\beta/(2\beta+1)}$

for all sufficiently large $n$.

**4. Adaptation over a general collection of convex parameter spaces.** Section 3 focused on collections of nested parameter spaces. It has been shown that the between class modulus of continuity completely characterizes the optimal expected length of adaptive confidence intervals. One particularly interesting feature of the nested case is that the optimal expected length of the confidence intervals does not depend on the number of parameter spaces in the collection.

The nested case, although interesting, is somewhat special. In this section general finite collections of convex parameter spaces are considered. In this general setting the theory is more complicated and in general the number of parameter spaces, say $k$, may also play a role in the optimal expected length of adaptive confidence intervals. For a fixed and finite number of parameter spaces the optimal expected length of adaptive intervals is still within a constant factor of the between class modulus of continuity. However the constant factor in this case can depend on the number of parameter spaces. We construct adaptive confidence intervals which show that this constant factor does not grow faster than $\sqrt{\log k}$ and we give an example which shows that this factor is sometimes necessary.

Let $\{\mathcal{F}_j : j = 1, \ldots, k\}$ be a collection of convex spaces with nonempty intersections, that is, $\mathcal{F}_i \cap \mathcal{F}_j \neq \varnothing$ for all $i$, $j$. The objective is to construct an adaptive confidence interval for a linear functional $Tf$ which has guaranteed coverage probability of $1 - \alpha$ over $\mathcal{G} = \bigcup_{j=1}^{k} \mathcal{F}_j$ and rate optimal expected length over each of the parameter spaces $\mathcal{F}_j$.

The adaptive confidence interval given in this section differs substantially from that given in the nested case. However, the general strategy for constructing adaptive confidence intervals in this setup is similar to that of the nested case. In particular, a key step is to first construct an interval which has optimal expected length over one of the parameter spaces while attaining coverage probability over the union of the parameter spaces.

4.1. *Constrained optimal expected length confidence intervals.* As mentioned above, it is convenient to construct a confidence interval which has shortest possible expected length over a given $\mathcal{F}_j$ while maintaining coverage probability over $\mathcal{G} = \bigcup_{j=1}^{k} \mathcal{F}_j$.

First note that for any confidence interval $CI \in \mathcal{I}_{\alpha, \mathcal{G}}$, Theorem 1 yields a target for the expected length

(53) $$L(CI, \mathcal{F}_j) \geq \left(\frac{1}{2} - \alpha\right) \omega_+ \left(\frac{z_\alpha}{\sqrt{n}}, \mathcal{F}_j, \mathcal{G}\right).$$

As in Section 2.2, for $1 \leq i, j \leq k$ set $\omega_{i,j} = \omega(\frac{z_{\alpha/2}}{\sqrt{n}}, \mathcal{F}_i, \mathcal{F}_j)$ and let $\hat{T}_{i,j}$ be a linear estimator which has variance bounded by $\frac{1}{z_{\alpha/2}^2}\omega_{i,j}^2$ and bias which satisfies

(54) $$\inf_{f \in \mathcal{F}_j}(E(\hat{T}_{i,j}) - Tf) \geq -\tfrac{1}{2}\omega_{i,j}$$

and

(55) $$\sup_{f \in \mathcal{F}_i}(E(\hat{T}_{i,j}) - Tf) \leq \tfrac{1}{2}\omega_{i,j}.$$

As a first step in the construction of adaptive confidence intervals, define $CI_{j,\alpha}^*$ by

(56) $$CI_{j,\alpha}^* = \left[\min_i\{\hat{T}_{i,j} - \tfrac{3}{2}\omega_{i,j}\}, \max_i\{\hat{T}_{j,i} + \tfrac{3}{2}\omega_{j,i}\}\right].$$

The following result shows that this confidence interval attains the lower bound on the maximum expected length over $\mathcal{F}_j$ given in (53) and satisfies the constraint that it has the minimum coverage of $1 - \alpha$ for all $f \in \mathcal{G}$.

PROPOSITION 4. *Let $\mathcal{F}_j$, $j = 1, \ldots, k$, be convex parameter spaces with $\mathcal{F}_i \cap \mathcal{F}_j \neq \varnothing$ for all $i, j$ and let $\mathcal{G} = \bigcup_{j=1}^k \mathcal{F}_j$. Let the interval $CI_{j,\alpha}^*$ be defined as in* (56). *Then $CI_{j,\alpha}^* \in \mathcal{I}_{\alpha,\mathcal{G}}$ and $CI_{j,\alpha}^*$ has expected length over $\mathcal{F}_j$ satisfying*

(57) $$L_\alpha^*(\mathcal{F}_j, \mathcal{G}) \leq L(CI_{j,\alpha}^*, \mathcal{F}_j) \leq \left\{\frac{8\sqrt{\log(k+1)} + 4z_{\alpha/2}}{(1/2 - \alpha)z_\alpha}\right\} \cdot L_\alpha^*(\mathcal{F}_j, \mathcal{G}).$$

REMARK 6. It follows from (59) that the expected length of the confidence interval $CI_{j,\alpha}^*$ is rate optimal as $n \to \infty$ as long as $k$ remains fixed.

PROOF OF PROPOSITION 4. First consider the coverage probability of the interval $CI_{j,\alpha}^*$. Suppose $f \in \mathcal{F}_m$ for some $1 \leq m \leq k$. Note that the interval $CI_{j,\alpha}^*$ contains

$$CI_{m,j} = [\hat{T}_{m,j} - \tfrac{3}{2}\omega_{m,j}, \hat{T}_{j,m} + \tfrac{3}{2}\omega_{j,m}].$$

The derivation below shows that the interval $CI_{m,j}$ has correct coverage probability. First note that for $f \in \mathcal{F}_m$, $E\hat{T}_{m,j} - Tf \leq \tfrac{1}{2}\omega_{m,j}$ and that $E\hat{T}_{j,m} - Tf \geq -\tfrac{1}{2}\omega_{j,m}$. Let

$$X_{m,j} = \frac{\hat{T}_{i,j} - Tf - (1/2)\omega_{m,j}}{\omega_{m,j}/z_{\alpha/2}},$$

$$X_{j,m} = \frac{\hat{T}_{j,m} - Tf + (1/2)\omega_{j,m}}{\omega_{j,m}/z_{\alpha/2}}.$$

Then for any $f \in \mathcal{F}_m$ it follows from (54) and (55) that $X_{m,j}$ has a Normal distribution with mean less than or equal to zero and variance bounded by 1 and $X_{j,m}$ has a Normal distribution with mean greater than or equal to zero and variance bounded by 1. Hence, for $f \in \mathcal{F}_m$,

$$\begin{aligned} P(Tf \in CI^*_{j,\alpha}) &\geq P(Tf \in CI_{m,j}) \\ &= P(X_{m,j} \geq -z_{\alpha/2} \text{ and } X_{j,m} \leq z_{\alpha/2}) \\ &\geq 1 - P(X_{m,j} \leq -z_{\alpha/2}) - P(X_{j,m} \geq z_{\alpha/2}) \\ &\geq 1 - \alpha. \end{aligned}$$

So for any $f \in \mathcal{G}$, $P(Tf \in CI^*_{j,\alpha}) \geq 1 - \alpha$ and thus coverage has been established.

The bounds on the expected length of these intervals can now be obtained by using the following technical lemma from Dudley [(1999), pages 56 and 57].

LEMMA 7. *Let $X_1, X_2, \ldots, X_k$ be normally distributed random variables with mean 0 and variance $\leq \sigma^2$. Then*

$$(58) \qquad E \max_{1 \leq i \leq k} |X_i| \leq \sigma \left(2 + \frac{4 + \log 4}{\log(3/2)}\right)^{1/2} \sqrt{\log(k+1)}.$$

Let

$$\xi_j = \omega_+\left(\frac{z_{\alpha/2}}{\sqrt{n}}, \mathcal{F}_j, \mathcal{G}\right) = \max_{1 \leq i \leq k} \{\omega_{i,j}, \omega_{j,i}\}.$$

It is easy to see that the length of the interval $CI^*_{j,\alpha}$ is bounded by

$$L(CI^*_{j,\alpha}) \leq \max_i (\hat{T}_{j,i} - Tf)_+ + \max_i (Tf - \hat{T}_{i,j})_+ + 3\xi_j.$$

Now note that if $f \in \mathcal{F}_j$, then for any $i \neq j$,

$$a_{j,i} \equiv E(\hat{T}_{j,i} - Tf) \leq \tfrac{1}{2}\omega_{j,i}$$

and

$$b_{i,j} \equiv E(Tf - \hat{T}_{i,j}) \leq \tfrac{1}{2}\omega_{i,j}.$$

Also note that for any real numbers $x$ and $y$, $(x+y)_+ \leq (x)_+ + (y)_+$. So for $f \in \mathcal{F}_j$ the expected length of $CI^*_{j,\alpha}$ satisfies

$$\begin{aligned} EL(CI^*_{j,\alpha}) &\leq E \max_i (\hat{T}_{j,i} - Tf)_+ + E \max_i (Tf - \hat{T}_{i,j})_+ + 3\xi_j \\ &\leq E\left(\max_i \{(a_{j,i})_+ + (\hat{T}_{j,i} - Tf - a_{j,i})_+\}\right) \end{aligned}$$

$$+ E\Big(\max_i\{(b_{i,j})_+ + (\hat{T}_{i,j} - Tf - b_{i,j})_+\}\Big) + 3\xi_j$$

$$\leq E\Big(\max_i(\hat{T}_{j,i} - \hat{T}f - a_{j,i})_+\Big) + E\Big(\max_i(\hat{T}_{i,j} - \hat{T}_j - b_{i,j})_+\Big) + 4\xi_j.$$

It then follows from Lemma 7 that

$$E_f(L(CI_{j,\alpha}^*)) \leq \frac{2}{z_{\alpha/2}}\xi_j\Big(2 + \frac{4 + \log 4}{\log(3/2)}\Big)^{1/2}\sqrt{\log(k+1)} + 4\xi_j$$

$$\leq \Big\{8\frac{\sqrt{\log(k+1)}}{z_{\alpha/2}} + 4\Big\}\xi_j$$

and it follows by taking the supremum over $\mathcal{F}_j$ that

(59) $$L(CI_{j,\alpha}^*, \mathcal{F}_j) \leq \Big\{8\frac{\sqrt{\log(k+1)}}{z_{\alpha/2}} + 4\Big\}\omega_+\Big(\frac{z_{\alpha/2}}{\sqrt{n}}, \mathcal{F}_j, \mathcal{G}\Big).$$

The proposition now follows by combining (20), (7) and (59). □

4.2. *Adaptive confidence intervals.* The intervals $CI_{j,\alpha}^*$ constructed in the last section have near optimal expected length over $\mathcal{F}_j$ but do not control the expected length over other $\mathcal{F}_i$. In this section adaptive confidence intervals over $\{\mathcal{F}_j : 1 \leq j \leq k\}$ are formed by intersecting such intervals. For a fixed $k$, the resulting interval has rate optimal expected length over every parameter space $\mathcal{F}_j$ for all $1 \leq j \leq k$. A Bonferroni approach is applied to the intervals of Section 4.1 to yield an adaptive confidence interval.

More specifically, define the confidence interval $CI^*$ by

(60) $$CI^* = \bigcap_{j=1}^k CI_{j,\alpha/k}^*$$

where $CI_{j,\alpha}^*$ are given in (56). The following theorem shows that this confidence interval has guaranteed coverage probability and also has near optimal expected length over $\mathcal{F}_j$ for each $1 \leq j \leq k$.

THEOREM 3. *Let $\mathcal{F}_j$, $j = 1, \ldots, k$, be convex parameter spaces with $\mathcal{F}_i \cap \mathcal{F}_j \neq \varnothing$ for all $i, j$ and let $\mathcal{G} = \bigcup_{j=1}^k \mathcal{F}_j$. Let the interval $CI^*$ be given as in (60). Then $CI^* \in \mathcal{I}_{\alpha,\mathcal{G}}$ and $CI^*$ satisfies*

(61) $$L_\alpha^*(\mathcal{F}_j, \mathcal{G}) \leq L(CI_\alpha^*, \mathcal{F}_j) \leq \frac{12 z_{\alpha/2k}}{(1/2 - \alpha)z_\alpha} \cdot L_\alpha^*(\mathcal{F}_j, \mathcal{G})$$

*for all $1 \leq j \leq k$.*

PROOF. The results follow easily from Proposition 4. For any $f \in \mathcal{G}$, Proposition 4 shows that
$$P(Tf \in CI^*_{j,\alpha/k}) \geq 1 - \frac{\alpha}{k}.$$

Hence, for any $f \in \mathcal{G}$,
$$P(Tf \in CI^*) = 1 - P(Tf \notin CI^*) \geq 1 - \sum_{j=1}^{k} P(Tf \notin CI^*_{j,\alpha}) \geq 1 - \alpha.$$

For the expected length note that
$$L(CI^*, \mathcal{F}_j) \leq L(CI^*_{j,\alpha/k}, \mathcal{F}_j) \leq \left\{ 8 \frac{\sqrt{\log(k+1)}}{z_{\alpha/2k}} + 4 \right\} \omega_+ \left( \frac{z_{\alpha/2k}}{\sqrt{n}}, \mathcal{F}_j, \mathcal{G} \right)$$

for any $1 \leq j \leq k$. For $0 < \alpha < 0.5$, calculations show that
$$\frac{\sqrt{\log(k+1)}}{z_{\alpha/2k}} \leq 1$$

and hence
$$(62) \qquad L(CI^*, \mathcal{F}_j) \leq 12 \omega_+ \left( \frac{z_{\alpha/2k}}{\sqrt{n}}, \mathcal{F}_j, \mathcal{G} \right).$$

The theorem now follows by combining (7), (20) and (62). □

REMARK 7. It follows from Lemma 6 that $z_{\alpha/2k} \leq \sqrt{\frac{2}{z_{\alpha/2}^2} \log k + 1} \cdot z_{\alpha/2}$. Hence it follows from (62) and (20) that
$$L(CI^*, \mathcal{F}_j) \leq 12 \omega_+ \left( \sqrt{\frac{2}{z_{\alpha/2}^2} \log k + 1} \cdot \frac{z_{\alpha/2}}{\sqrt{n}}, \mathcal{F}_j, \mathcal{G} \right)$$
$$\leq 12 \sqrt{\frac{2}{z_{\alpha/2}^2} \log k + 1} \cdot \omega_+ \left( \frac{z_{\alpha/2}}{\sqrt{n}}, \mathcal{F}_j, \mathcal{G} \right).$$

The ratio of the upper bound just given to the lower bound in (53) is thus clearly bounded by a constant multiple of $\sqrt{\log k}$.

Section 5.2 gives an example of a nearly black object which shows that this $\sqrt{\log k}$ factor cannot in general be improved.

**5. Minimax confidence interval for nonconvex parameter spaces.** As mentioned in the Introduction, Donoho (1994) constructed for any convex parameter space $\mathcal{F}$ fixed length intervals centered at affine estimators which have length within a small constant factor of the minimax expected length $L_\alpha^*(\mathcal{F})$. Although the focus of the present paper is on adaptation the adaptation theory developed in the previous sections can also be used to yield a minimax theory for parameter spaces that are finite unions of convex parameter spaces. In this section confidence intervals with a specified coverage probability are given which also have near optimal maximum expected length. It is also shown, in contrast to the theory for convex parameter spaces, that optimal confidence intervals centered on affine estimators can have expected length much longer than the expected length of optimal confidence intervals centered at nonlinear estimators.

Let $\mathcal{F}_i$, $i = 1, \ldots, k$, be convex parameter spaces with $\mathcal{F}_i \cap \mathcal{F}_j \neq \varnothing$ for all $i$, $j$ and let $\mathcal{G} = \bigcup_{i=1}^k \mathcal{F}_i$. Note that the parameter space $\mathcal{G}$ is in general nonconvex. The minimax expected length of confidence intervals $CI \in \mathcal{I}_{\alpha,\mathcal{G}}$ can be bounded above and below as follows.

Set $0 < \alpha < \frac{1}{2}$ and let $CI$ be a $1 - \alpha$ level confidence interval for all $f \in \mathcal{G} = \bigcup_{i=1}^k \mathcal{F}_i$. It follows from Theorem 1 that the maximum expected length of $CI \in \mathcal{I}_{\alpha,\mathcal{G}}$ is bounded below by

$$(63) \qquad L(CI, \mathcal{G}) \geq \left(\frac{1}{2} - \alpha\right) \omega\left(\frac{z_\alpha}{\sqrt{n}}, \mathcal{G}\right).$$

Upper bounds on the minimax expected length can be obtained by considering the confidence interval $CI^*$ as defined in (60). As shown in Theorem 3 this interval has coverage probability of at least $1 - \alpha$ over $\mathcal{G}$. In addition, it follows from (61) that the maximum of the expected length of $CI^*$ over $\mathcal{G}$ satisfies

$$(64) \qquad \begin{aligned} L(CI^*, \mathcal{G}) &= \max_{1 \leq j \leq k} L(CI^*, \mathcal{F}_j) \\ &\leq 12 \max_{1 \leq j \leq k} \omega_+\left(\frac{z_{\alpha/2k}}{\sqrt{n}}, \mathcal{F}_j, \mathcal{G}\right) \\ &= 12 \omega\left(\frac{z_{\alpha/2k}}{\sqrt{n}}, \mathcal{G}\right). \end{aligned}$$

Hence, (63) and (64) together yield the following result on the minimax expected length of $1 - \alpha$ level confidence intervals over $\mathcal{G}$.

THEOREM 4. *Let $\mathcal{G} = \bigcup_{j=1}^k \mathcal{F}_i$, where for $j = 1, \ldots, k$ $\mathcal{F}_j$ are convex spaces with nonempty intersections and suppose $0 < \alpha < \frac{1}{2}$. Then*

$$(65) \qquad \left(\frac{1}{2} - \alpha\right) \omega\left(\frac{z_\alpha}{\sqrt{n}}, \mathcal{G}\right) \leq L_\alpha^*(\mathcal{G}) \leq 12 \omega\left(\frac{z_{\alpha/2k}}{\sqrt{n}}, \mathcal{G}\right).$$

Hence, the confidence interval $CI^*$ attains the optimal rate of convergence for the maximum expected length over the parameter spaces $\mathcal{G}$ when the number of convex subspaces is fixed and finite.

The example of confidence intervals in Section 5.2 for a linear functional of nearly black objects shows that the factor of $z_{\alpha/2k} \asymp \sqrt{\log k}$ in the upper bound of (65) cannot be dropped in general when the number $k$ of convex subspaces grows with $n$.

5.1. *Confidence intervals centered at affine estimators.* We now consider the performance of confidence intervals centered at affine estimators over nonconvex parameter spaces. As mentioned earlier, when the parameter space $\mathcal{F}$ is assumed to be fixed and convex, Donoho (1994) and Theorem 1 given in Section 2 together show that the length of the shortest fixed length confidence interval centered on an affine estimator is within a fixed constant factor of the maximum expected length of the optimal confidence interval. Hence there is relatively little to gain by looking beyond the class of fixed length confidence intervals centered on affine estimators.

The following theorem considers the case when the parameter space is nonconvex. Once again let C.Hull($\mathcal{F}$) denote the convex hull of a parameter space $\mathcal{F}$.

THEOREM 5. *Consider the white noise model* (1) *or the sequence model* (2). *Let $\hat{T}$ be an affine estimator of $Tf$ and $\gamma \geq 0$ a nonnegative random variable. If $CI = [\hat{T} - \gamma, \hat{T} + \gamma]$ is a (variable length) confidence interval centered at $\hat{T}$ and $CI \in \mathcal{I}_{\alpha,\mathcal{F}}$, then*

$$(66) \qquad L(CI,\mathcal{F}) \geq C(\alpha)\omega\left(\frac{2z_{\alpha/2}}{\sqrt{n}}, \text{C.Hull}(\mathcal{F})\right)$$

*where $C(\alpha) > 0$ is a constant depending on $\alpha$ only. In particular, if the interval $CI$ is of fixed length, then*

$$(67) \qquad L(CI) \geq \frac{1}{2}\omega\left(\frac{2z_{\alpha/2}}{\sqrt{n}}, \text{C.Hull}(\mathcal{F})\right).$$

PROOF. It is shown in Cai and Low (2004) that the affine estimator $\hat{T}$ satisfies

$$\sup_{f \in \mathcal{F}} |E\hat{T} - Tf| = \sup_{f \in \text{C.Hull}(\mathcal{F})} |E\hat{T} - Tf|.$$

It then follows from Theorem 2 of Low (1995) that $\hat{T}$ must satisfy either

$$(68) \qquad \sup_{f \in \mathcal{F}} |E\hat{T} - Tf| \geq \frac{1}{4}\omega\left(\frac{2z_{\alpha/2}}{\sqrt{n}}, \text{C.Hull}(\mathcal{F})\right)$$

or

$$(69) \qquad \sigma_{\hat{T}} \geq \frac{1}{4z_{\alpha/2}}\omega\left(\frac{2z_{\alpha/2}}{\sqrt{n}}, \text{C.Hull}(\mathcal{F})\right),$$

where $\sigma_{\hat{T}}$ denotes the standard deviation of the estimator $\hat{T}$. We now consider the two cases separately. If (68) holds, then for any $\varepsilon > 0$, there exists $f \in \mathcal{F}$ such that

$$(70) \qquad B_f \equiv |E\hat{T} - Tf| \geq \frac{1}{4}\omega\left(\frac{2z_{\alpha/2}}{\sqrt{n}}, \text{C.Hull}(\mathcal{F})\right) - \varepsilon.$$

Since $CI = [\hat{T} - \gamma, \hat{T} + \gamma]$ has minimum coverage probability of at least $1 - \alpha$ over $\mathcal{F}$,

$$1 - \alpha \leq P_f(|\hat{T} - Tf| \leq \gamma)$$
$$= P_f(|\hat{T} - Tf| \leq \gamma \text{ and } \gamma \leq B_f) + P_f(|\hat{T} - Tf| \leq \gamma \text{ and } \gamma > B_f)$$
$$\leq P_f(|\hat{T} - Tf| \leq B_f) + P(\gamma > B_f).$$

Since $\hat{T}$ is an affine estimator and thus has a normal distribution, it is easy to check that $P_f(|\hat{T} - Tf| \leq B_f) \leq 1/2$ and hence

$$(71) \qquad P(\gamma > B_f) \geq \tfrac{1}{2} - \alpha.$$

Letting $\varepsilon \to 0$ in (70), it then follows that

$$E_f L(CI) = 2E_f(\gamma) \geq 2B_f P(\gamma > B_f)$$
$$\geq \left(\frac{1}{4} - \frac{\alpha}{2}\right)\omega\left(\frac{2z_{\alpha/2}}{\sqrt{n}}, \text{C.Hull}(\mathcal{F})\right).$$

If (69) holds, we have, for $f \in \mathcal{F}$,

$$1 - \alpha \leq P_f(|\hat{T} - Tf| \leq \gamma)$$
$$= P\left(-\frac{\gamma}{\sigma_{\hat{T}}} - \frac{E\hat{T} - Tf}{\sigma_{\hat{T}}} \leq Z \leq \frac{\gamma}{\sigma_{\hat{T}}} - \frac{E\hat{T} - Tf}{\sigma_{\hat{T}}}\right)$$
$$\leq P\left(|Z| \leq \frac{\gamma}{\sigma_{\hat{T}}}\right)$$
$$= P\left(|Z| \leq \frac{\gamma}{\sigma_{\hat{T}}} \text{ and } \gamma \leq z_{0.25}\sigma_{\hat{T}}\right) + P\left(|Z| \leq \frac{\gamma}{\sigma_{\hat{T}}} \text{ and } \gamma > z_{0.25}\sigma_{\hat{T}}\right)$$
$$\leq P(|Z| \leq z_{0.25}) + P(\gamma > z_{0.25}\sigma_{\hat{T}})$$

where $Z$ denotes a standard normal random variable. Hence

$$P(\gamma > z_{0.25}\sigma_{\hat{T}}) \geq \tfrac{1}{2} - \alpha.$$

Consequently,

$$E_f L(CI) = 2E_f(\gamma) \geq 2z_{0.25}\sigma_{\hat{T}} P(\gamma > z_{0.25}\sigma_{\hat{T}})$$
$$\geq (1 - 2\alpha)\frac{z_{0.25}}{4z_{\alpha/2}}\omega\left(\frac{2z_{\alpha/2}}{\sqrt{n}}, \text{C.Hull}(\mathcal{F})\right)$$

$$\geq \frac{1-2\alpha}{10z_{\alpha/2}}\omega\bigg(\frac{2z_{\alpha/2}}{\sqrt{n}}, \text{C.Hull}(\mathcal{F})\bigg).$$

Equation (66) now follows by taking $C(\alpha) = \min\{\frac{1}{4} - \frac{\alpha}{2}, \frac{1-2\alpha}{10z_{\alpha/2}}\}$. Equation (67) for the fixed length case is easier to prove and we omit the proof here. □

REMARK 8. Theorem 5 shows that the minimax expected length of confidence intervals centered at affine estimators is determined by the modulus of continuity over the convex hull of $\mathcal{F}$, not over $\mathcal{F}$ itself. In the case that $\omega(\varepsilon, \text{C.Hull}(\mathcal{F})) \gg \omega(\varepsilon, \mathcal{F})$, any confidence intervals centered at affine estimators will perform poorly. Such is the case in the near black object example given in the next section.

5.2. *Nearly black object.* In this section an example is given which shows that the factor $z_{\alpha/2k} \asymp \sqrt{\log k}$ in the upper bound of the minimax expected length given in Theorem 4 cannot in general be dropped. It is also shown that confidence intervals centered at affine estimators are far from optimal.

Consider the Gaussian sequence model (2) with the index set $\mathcal{M} = \{1, 2, \ldots, n\}$, namely

(72) $$Y(i) = f(i) + n^{-1/2}z_i, \qquad i = 1, \ldots, n,$$

where $z_i \overset{\text{i.i.d.}}{\sim} N(0,1)$. The size of the vector, $n$, is assumed large. We assume that the vector $f$ is sparse: only a small fraction of components are nonzero, and the indices or locations of the nonzero components are not known in advance.

Denote the $\ell_0$ quasi-norm by $\|f\|_0 = \text{Card}(\{i: f(i) \neq 0\})$. Fix $m_n$. The collection of vectors with at most $m_n$ nonzero entries is

$$\mathcal{G} = \ell_0(m_n) = \{f \in \mathbb{R}^n : \|f\|_0 \leq m_n\}.$$

Assume that $m_n$ is known and $m_n \leq n^\gamma$ where $\gamma < \frac{1}{2}$.

Such an example is considered in Cai and Low (2004) in the context of minimax estimation. The model, which arises naturally in wavelet analysis, has also been studied in Donoho, Johnstone, Hoch and Stern (1992) and Abramovich, Benjamini, Donoho and Johnstone (2000) for estimating the whole object.

Let the linear functional $Tf$ be given by

$$Tf = \sum_{i=1}^{n} f(i),$$

and following Cai and Low (2004) let $\mathcal{I}(m_n, n)$ be the class of all subsets of $\{1, \ldots, n\}$ of $m_n$ elements and for $I \in \mathcal{I}(m_n, n)$ let

$$\mathcal{F}_I = \{f \in \mathbb{R}^n : f(j) = 0 \ \forall j \notin I\}.$$

Note that $\mathcal{F}_I$ is an $m_n$-dimensional subspace spanned by the coordinates in $I$. These are obviously convex and $\mathcal{G} = \cup \mathcal{F}_I$ where the union is taken over $I$ in the set $\mathcal{I}(m_n, n)$. From now on we shall assume that $I$ is in the set $\mathcal{I}(m_n, n)$.

Simple calculations show that for all $I, J \in \mathcal{I}(m_n, n)$

$$\omega(\varepsilon, \mathcal{F}_I, \mathcal{F}_J) = \sqrt{\text{Card}(I \cup J)}\varepsilon$$

and consequently

$$\omega(\varepsilon, \mathcal{F}_I, \mathcal{G}) = \omega(\varepsilon, \mathcal{G}, \mathcal{F}_I) = \omega(\varepsilon, \mathcal{G}) = \sqrt{2m_n}\varepsilon.$$

REMARK 9. It is easy to see that $\text{C.Hull}(\mathcal{G}) = \mathbb{R}^n$ and hence $\omega(\varepsilon, \text{C.Hull}(\mathcal{G})) = \sqrt{n}\varepsilon$. It follows from Theorem 5 and (66) that any confidence interval with coverage of at least $1 - \alpha$ centered at an affine estimator must have maximum expected length bounded from below by a fixed constant not depending on $n$.

Let $k$ be the number of the $m_n$-dimensional parameter spaces $\mathcal{F}_I$. Then $k$ is equal to $n$ choose $m_n$ and it is easy to see that

$$k = \binom{n}{m_n} \leq n^{m_n}.$$

The following result gives a lower bound on the expected length of any confidence interval with a minimum coverage probability of $1 - \alpha$ over $\mathcal{G}$.

PROPOSITION 5. *Suppose that we observe the Gaussian sequence model* (72), *that* $n \geq 4$ *and* $m_n < n^\gamma$ *with* $\gamma < \frac{1}{2}$. *Let* $Tf = \sum_{i=1}^n f(i)$ *and* $0 < \alpha < \frac{1}{2}$. *Suppose that* $CI(Y)$ *is a confidence interval for* $Tf$ *based on* (72) *and* $CI(Y) \in \mathcal{I}_{\alpha, \mathcal{G}}$. *Then for all sufficiently large* $n$,

$$\begin{aligned}
E_0 L(CI(Y)) &\geq \left(\frac{1}{2} - \alpha\right) \frac{m_n}{\sqrt{n}} \sqrt{\frac{1}{2} \log\left(\frac{n}{m_n^2}\right)} \\
&\geq \left(\frac{1}{2} - \alpha\right) \sqrt{\frac{1}{4} - \frac{\gamma}{2}} \cdot \omega\left(\frac{\sqrt{\log k}}{\sqrt{n}}, \mathcal{G}\right)
\end{aligned}$$
(73)

*where* $E_0$ *denotes expectation under the Gaussian model* (72) *with* $f(i) = 0$ *for* $i = 1, 2, \ldots, n$.

REMARK 10. It follows immediately from (73) that the maximum expected length of $CI(Y)$ over $\mathcal{G}$ satisfies

$$L(CI(Y), \mathcal{G}) \geq C\omega\left(\frac{\sqrt{\log k}}{\sqrt{n}}, \mathcal{G}\right).$$
(74)

Comparing the lower bound (74) for the maximum expected length with the minimax upper bound given in (65) shows that the factor $\sqrt{\log k}$ in the upper bound for the minimax expected length cannot be dropped in general. A similar result also holds for adaptation.

PROOF OF PROPOSITION 5. In the following proof the calculation of the $L_1$ distance between a mixture of normals and a given normal distribution follows a similar calculation used in Cai and Low (2004). We include the details of the calculation here for completeness. In the proof we will omit the subscript in $m_n$ and simply write $m$ for $m_n$. Let $\psi_f$ be the joint density of the Gaussian observations given in (72). More specifically $\psi_f$ is a multivariate normal density with mean $(f(1), f(2), \ldots, f(n))$ and covariance matrix $\frac{1}{n}A_n$ where $A_n$ is the $n \times n$ identity matrix. Fix a constant $\rho > 0$. For $I \in \mathcal{I}(m,n)$ let $f_I$ be defined by $f_I(j) = \frac{\rho}{\sqrt{n}}\mathbb{1}(j \in I)$ and let $f_0$ be the sequence defined by $f_0(j) = 0$ for $j = 1, 2, \ldots, n$. Finally let

$$\psi_* = \frac{1}{\binom{n}{m}} \sum_{I \in \mathcal{I}(m,n)} \psi_{f_I}.$$

Note that a similar mixture prior was used in Baraud (2002) to give lower bounds in a nonparametric testing problem. Note that for all $f_I$, $Tf_I = m\frac{\rho}{\sqrt{n}}$ and that $Tf_0 = 0$. Note also that if

$$P_{\psi_{f_I}}\left(m\frac{\rho}{\sqrt{n}} \in CI(Y)\right) \geq 1 - \alpha$$

for all $I \in \mathcal{I}(m,n)$ then it follows that

$$P_{\psi_*}\left(m\frac{\rho}{\sqrt{n}} \in CI(Y)\right) = \frac{1}{\binom{n}{m}} \sum_{I \in \mathcal{I}(m,n)} \psi_{f_I} P_{\psi_{f_I}}\left(m\frac{\rho}{\sqrt{n}} \in CI(Y)\right) \geq 1 - \alpha.$$

Note that

$$\int \frac{\psi_*^2}{\psi_{f_0}} = \frac{1}{\binom{n}{m}^2} \sum_{I \in \mathcal{I}(m,n)} \sum_{I' \in \mathcal{I}(m,n)} \int \frac{\psi_{f_I}\psi_{f_{I'}}}{\psi_{f_0}},$$

and simple calculations show that

$$\int \frac{\psi_{f_I}\psi_{f_{I'}}}{\psi_{f_0}} = \exp(j\rho^2),$$

where $j$ is the number of points in the set $I \cap I'$. It follows that

$$\int \frac{\psi_*^2}{\psi_{f_0}} = E\exp(J\rho^2),$$

where $J$ has a hypergeometric distribution,

$$P(J=j) = \frac{\binom{m}{j}\binom{n-m}{m-j}}{\binom{n}{m}}.$$

Now note that from Feller [(1968), page 59],

$$P(J=j) \leq \binom{m}{j}\left(\frac{m}{n}\right)^j\left(1-\frac{m}{n}\right)^{m-j}\left(1-\frac{m}{n}\right)^{-m}.$$

Now suppose that $n \geq 4$ and that $m < n^{1/2}$. Then

$$\left(1-\frac{m}{n}\right)^{-m} \leq 4^{m^2/n}$$

and hence

$$P(J=j) \leq 4^{m^2/n}\binom{m}{j}\left(\frac{m}{n}\right)^j\left(1-\frac{m}{n}\right)^{m-j}.$$

It now follows that if $n \geq 4$ and $m < n^\gamma$ with $\gamma < \frac{1}{2}$, then

$$\int \frac{\psi_*^2}{\psi_{f_0}} = E\exp(J\rho^2)$$

$$\leq 4^{m^2/n}\left(1 - \frac{m}{n} + \frac{m}{n}\exp(\rho^2)\right)^m$$

$$\leq 4^{m^2/n}\left(1 + \frac{m}{n}\exp(\rho^2)\right)^m.$$

Now take $\rho = \sqrt{\frac{1}{2}\log\frac{n}{m^2}}$. Then

$$\int \frac{\psi_*^2}{\psi_{f_0}} \leq 4^{n^{-(1-2\gamma)}}\left(1 + \frac{1}{n^{1/2}}\right)^{n^\gamma} \downarrow 1.$$

Hence we can bound the $L_1$ distance by

$$\int |\psi_* - \psi_{f_0}| \leq \left(\int \frac{(\psi_* - \psi_{f_0})^2}{\psi_{f_0}}\right)^{1/2} = \left(\int \frac{\psi_*^2}{\psi_{f_0}} - 1\right)^{1/2}$$

$$\leq \left(4^{n^{-(1-2\gamma)}}\left(1 + \frac{1}{n^{1/2}}\right)^{n^\gamma} - 1\right)^{1/2} \downarrow 0.$$

So for any $0 < \varepsilon < 1 - 2\alpha$ there exists $n_\varepsilon$ such that for all $n \geq n_\varepsilon$, $\int |\psi_* - \psi_{f_0}| \leq \varepsilon$.

It follows from the fact that $CI$ has minimum coverage probability of $1 - \alpha$ and that the $L_1$ distance between $\psi_{f_0}$ and $\psi_*$ is bounded above by $\varepsilon$ that

$$P_{\psi_{f_0}}\left(m\frac{\rho}{\sqrt{n}} \in CI\right) \geq P_{\psi_*}\left(m\frac{\rho}{\sqrt{n}} \in CI\right) - \varepsilon \geq 1 - \alpha - \varepsilon.$$

Hence

$$P_{\psi_{f_0}}\left(0 \in CI \text{ and } m\frac{\rho}{\sqrt{n}} \in CI\right) \geq 1 - 2\alpha - \varepsilon.$$

Since $CI$ is an interval the length of this interval must be at least $m\frac{\rho}{\sqrt{n}}$ when both 0 and $m\frac{\rho}{\sqrt{n}}$ are in the interval. Hence for $n \geq n_\varepsilon$,

$$E_{\psi_{f_0}} L(CI(Y)) \geq (1 - 2\alpha - \varepsilon)\frac{m}{\sqrt{n}}\sqrt{\frac{1}{2}\log\left(\frac{n}{m^2}\right)}.$$

Now take $\varepsilon = \frac{1}{2} - \alpha$. Then for all sufficiently large $n$,

$$\begin{aligned}E_{\psi_{f_0}} L(CI(Y)) &\geq \left(\frac{1}{2} - \alpha\right)\frac{m}{\sqrt{n}}\sqrt{\frac{1}{2}\log\left(\frac{n}{m^2}\right)} \\ &\geq \left(\frac{1}{2} - \alpha\right)\sqrt{\frac{1}{4} - \frac{\gamma}{2}} \cdot \omega\left(\frac{\sqrt{\log k}}{\sqrt{n}}, \mathcal{G}\right),\end{aligned}$$

where $k$ is the number of convex parameter spaces in $\mathcal{G}$. $\square$

REMARK 11. It follows immediately from Proposition 5 that

$$L_\alpha^*(\mathcal{F}_I) \geq \left(\frac{1}{2} - \alpha\right)\sqrt{\frac{1}{4} - \frac{\gamma}{2}} \cdot \omega_+\left(\frac{\sqrt{\log k}}{\sqrt{n}}, \mathcal{F}_I, \mathcal{G}\right).$$

Hence the factor of $z_{\alpha/2k} \asymp \sqrt{\log k}$ for adaptation in the upper bound of Theorem 3 and the same factor for minimax confidence procedures in Theorem 4 cannot in general be removed.

**Acknowledgment.** We thank one of the referees for very thorough and useful comments which have helped to improve the presentation of the paper. In particular, some of the notation used in the paper is based on suggestions made by this referee.

DEPARTMENT OF STATISTICS
THE WHARTON SCHOOL
UNIVERSITY OF PENNSYLVANIA
PHILADELPHIA, PENNSYLVANIA 19104
USA
E-MAIL: tcai@wharton.upenn.edu